\documentclass[a4paper]{amsart}
\usepackage{amssymb}
\usepackage{amscd}
\usepackage{amsthm}
\usepackage{amsmath}
\usepackage{mathtools}
\usepackage{latexsym}
\usepackage{hyperref}
\usepackage[all]{xy}
\usepackage[utf8]{inputenc}
\usepackage{enumitem}

\usepackage{tikz-cd}

\setlist[enumerate]{label={(\roman*)}}

\theoremstyle{plain}
\newtheorem{theorem}{Theorem}
\newtheorem{corollary}[theorem]{Corollary}
\newtheorem{lemma}[theorem]{Lemma}
\newtheorem{proposition}[theorem]{Proposition}

\theoremstyle{definition}
\newtheorem{definition}[theorem]{Definition}
\newtheorem{example}[theorem]{Example}

\theoremstyle{remark}

\newtheorem{remark}{Remark}

\newtheorem{conjecture}[theorem]{Conjecture}
\numberwithin{theorem}{section}
\usepackage{newtxtext}
\usepackage{newtxmath}

\DeclareMathAlphabet\urwscr{U}{urwchancal}{m}{n}%
\DeclareMathAlphabet\rsfscr{U}{rsfso}{m}{n}
\DeclareMathAlphabet\euscr{U}{eus}{m}{n}
\DeclareFontEncoding{LS2}{}{}
\DeclareFontSubstitution{LS2}{stix}{m}{n}
\DeclareMathAlphabet\stixcal{LS2}{stixcal}{m} {n}
\newcommand{\Filt}[1]{\mbox{\rm Filt}(#1)}
\newcommand{\Qcoh}[1]{\mbox{\rm Qcoh}(#1)}
\newcommand{\Spec}[1]{\operatorname{Spec}(#1)}

\newcommand{\card}{\mbox{\rm{card\,}}}

\newcommand{\Add}{\mbox{\rm{Add\,}}}

\newcommand{\Def}{\mbox{\rm{Def\,}}}

\newcommand{\End}{\mbox{\rm{End\,}}}

\newcommand{\im}{\mbox{\rm{Im\,}}}
\newcommand{\Ker}{\mbox{\rm{Ker\,}}}

\newcommand{\Hom}[3]{\operatorname{Hom}_{#1}(#2,#3)}
\newcommand{\Ext}[4]{\operatorname{Ext}^{#1}_{#2}(#3,#4)}
\newcommand{\Tor}[4]{\mbox{\rm{Tor}}_{#1}^{#2}(#3,#4)}

\newcommand{\rmod}[1]{\mbox{\rm{Mod}--}{#1}}

\newcommand{\lmod}[1]{{#1}\mbox{--\rm{Mod}}}

\begin{document}

\title{Flat Mittag-Leffler modules, and their relative and restricted versions}

\author{Jan Trlifaj}
\address{Charles University, Faculty of Mathematics
and Physics, Department of Algebra \\
Sokolovsk\'{a} 83, 186 75 Prague 8, Czech Republic}
\email{trlifaj@karlin.mff.cuni.cz}

\begin{abstract} Assume that $R$ is a non-right perfect ring. Then there is a proper class of classes of (right $R$-) modules closed under transfinite extensions lying between the classes $\mathcal P _0$ of projective modules, and $\mathcal F _0$ of flat modules. These classes can be defined as variants of the class $\mathcal F \mathcal M$ of absolute flat Mittag-Leffler modules: either as their restricted versions (lying between $\mathcal P _0$ and $\mathcal F \mathcal M$), or their relative versions (between $\mathcal F \mathcal M$ and $\mathcal F _0$). In this survey, we will deal with applications of these classes in relative homological algebra and algebraic geometry.

The classes $\mathcal P _0$ and $\mathcal F _0$ are known to provide for approximations, and minimal approximations, respectively. We will show that the classes of restricted flat relative Mittag-Leffler modules, and flat relative Mittag-Leffler modules, have rather different approximation properties: the former classes always provide for approximations, but the latter do not, except for the boundary case of $\mathcal F _0$. 

The notion of an (infinite dimensional) vector bundle is known to be Zariski local for all schemes, the key point of the proof being that projectivity ascends and descends along flat and faithfully flat ring homomorphisms, respectively. We will see that the same holds for the properties of being a $\kappa$-restricted flat Mittag-Leffler module for each cardinal $\kappa \geq \aleph_0$, and also a flat $\mathcal Q$-Mittag-Leffler module whenever $\mathcal Q$ is a definable class of finite type. Thus, as in the model case of vector bundles, Zariski locality holds for flat quasi-coherent sheaves induced by each of these classes of modules. Moreover, we will see that the notion of a locally $n$-tilting quasi-coherent sheaf is Zariski local for all $n \geq 0$. 
\end{abstract}

\date{\today}

\thanks{Research supported by GA\v CR 23-05148S}
	

\maketitle

\section{Introduction}

By a classic result of Bass, non-right perfect rings $R$ are characterized by the existence of countably presented flat (right $R$-) modules that are not projective \cite[28.4]{AF}. While projective modules can always be decomposed into direct sums of countably generated submodules \cite[26.2]{AF}, only a weak decomposition theorem is available for the flat modules: if $\kappa = \card R + \aleph_0$, then each flat module $M$ can be deconstructed into a \emph{transfinite extension} of $\leq \kappa$-presented flat modules \cite[6.17]{GT}. That is, $M$ possesses a continuous increasing chain of submodules, $( M_\alpha \mid \alpha \leq \sigma)$, such that $M_0 = 0$, $M_\sigma = M$, and for each $\alpha < \sigma$, $M_{\alpha + 1}/M_\alpha$ is a $\leq \kappa$-presented flat module. 

Motivated by Grothendieck's question on Zariski locality of the notion of a vector bundle, Raynaud and Gruson introduced the intermediate class of (absolute) flat Mittag-Leffler modules, $\mathcal F \mathcal M$, in \cite{RG}. Recall that a module $M$ is \emph{Mittag-Leffler}, if for each family $\mathcal I = ( Q_i \mid i \in I )$ of left $R$-modules, the canonical group homomorphism $\varphi _{\mathcal I} : M \otimes_R \prod_{i \in I} Q_i \to \prod_{i \in I} M \otimes_R  Q_i$ is monic (see below for unexplained terminology). 

If $R$ is not right perfect, then $\mathcal P_0 \subsetneq \mathcal F \mathcal M \subsetneq \mathcal F _0$ where $\mathcal P_0$ and $\mathcal F _0$ denote the classes of all projective and flat modules, respectively. All these classes are closed under transfinite extensions, but unlike $\mathcal P_0$ and $\mathcal F_0$, the class $\mathcal F \mathcal M$ is not deconstructible \cite[10.13]{GT}. However, there is a very rich supply of deconstructible classes closed under transfinite extensions between $\mathcal P_0$ and $\mathcal F_0$: as $\kappa$ ranges over all infinite cardinals, the classes $\mathcal F \mathcal M _\kappa$ of \emph{$\kappa$-restricted flat Mittag-Leffler modules} (= transfinite extensions of $\leq \kappa$-presented flat Mittag-Leffler modules) form a strictly increasing chain $( \mathcal F \mathcal M _\kappa \mid \aleph_0 \leq \kappa )$ between $\mathcal P_0$ and $\mathcal F \mathcal M$:
$$\mathcal P _0 = \mathcal F \mathcal M _{\aleph_0} \subsetneq \mathcal F \mathcal M _{\aleph_1} \subsetneq \dots \subsetneq \mathcal F \mathcal M _\kappa \subsetneq \mathcal F \mathcal M _{\kappa^+} \subsetneq \dots \subsetneq \bigcup_{\aleph_0 \leq \kappa} \mathcal F \mathcal M _\kappa = \mathcal F \mathcal M.$$                  
For each $\kappa \geq \aleph_0$, the class $\mathcal F \mathcal M _\kappa$ is obviously deconstructible, and hence precovering \cite[7.21]{GT}, but the class $\mathcal F \mathcal M$ fails these properties \cite[3.3]{S}. 

When $R$ is not right perfect, there is also a rich intermediate structure between the classes $\mathcal F \mathcal M$ and $\mathcal F _0$ provided by the classes of flat \emph{relative} Mittag-Leffler modules. These are obtained by restricting the choice of the families $\mathcal I$ in the definition above: if $\mathcal Q$ is any class of left $R$-modules, then a module $M$ is \emph{$\mathcal Q$-Mittag-Leffler}, if the canonical group homomorphism $\varphi _{\mathcal I}$ is monic for each family $\mathcal I = ( Q_i \mid i \in I )$ which consists of modules from $\mathcal Q$. Following \cite{HT}, we will denote the class of all flat $\mathcal Q$-Mittag-Leffler modules by $\mathcal D _{\mathcal Q}$. So if $\mathcal Q ^\prime \subseteq \mathcal Q$, we have the inclusions
$$\mathcal F \mathcal M = \mathcal D _{\lmod R} \subseteq \mathcal D _{\mathcal Q} \subseteq \mathcal D _{\mathcal Q ^\prime} \subseteq \mathcal F _0.$$   

Though there is a proper class of classes $\mathcal Q \subseteq \lmod R$, there is only a set of different classes $\mathcal D _{\mathcal Q}$. As proved by Rothmaler \cite[2.2]{R1}, $\mathcal D _{\mathcal Q} = \mathcal D _{\Def {(\mathcal Q)}}$ where $\Def {(\mathcal Q)}$ is the \emph{definable closure} of $\mathcal Q$, that is, the least class of left $R$-modules containing $\mathcal Q$ and closed under direct products, direct limits, and pure submodules. Moreover, if $R \in \mathcal Q $, then the structure of the class $\mathcal D _{\mathcal Q}$ is completely determined by the countably presented modules in $\mathcal D _{\mathcal Q}$. So if $R \in \mathcal Q ^\prime$, then $\mathcal D _{\mathcal Q} = \mathcal D _{\mathcal Q ^\prime}$, iff $\mathcal D _{\mathcal Q}$ and $\mathcal D _{\mathcal Q ^\prime}$ contain the same countably presented modules, \cite[2.5]{BT1}. As for the approximation properties of flat relative Mittag-Leffler modules, the situation is similar to the absolute case: the class $\mathcal D _{\mathcal Q}$ is precovering only if it coincides with the class of all flat modules \cite[2.6]{BT1}.  

Originally, absolute flat Mittag-Leffler modules served as a tool for proving Zariski locality of the notion of a vector bundle over any scheme, cf.\ \cite{P} and \cite{RG}. Relative Mittag-Leffler modules turned out to play an important role in (infinite dimensional) tilting theory, \cite{AH}. This has led to an investigation of quasi-coherent sheaves associated with tilting. Their Zariski locality was proved in \cite{HST}. The Zariski locality for quasi-coherent sheaves induced by restricted flat Mittag-Leffler modules goes back to \cite{EGT}, while the corresponding result for quasi-coherent sheaves induced by (some) flat relative Mittag-Leffler modules is quite recent, \cite{BT2}.  

\medskip
The goal of this survey is to provide a unified presentation of these recent results: in Section \ref{approx}, we will deal with the approximation properties, while Section \ref{Zariski} concerns the Zariski locality of the various induced notions of quasi-coherent sheaves on schemes.             

\section{Basic notions and examples}     

Mittag-Leffler modules are closely related to Mittag-Leffler inverse systems of modules. Thus we start with recalling basic notions, and fixing our notation, concerning direct and inverse limits of direct and inverse systems of modules over an arbitrary ring $R$. We will also take the opportunity to present examples showing the variety of properties of these systems. 

For a class $\mathcal C \subseteq \rmod R$, we will denote by $\Filt {\mathcal C}$ the class of all modules that are \emph{transfinite extensions} of the modules from $\mathcal C$ (or \emph{$\mathcal C$-filtered modules}), that is, the modules $M$ possessing an increasing chain of submodules, $\mathcal M = ( M_\alpha \mid \alpha \leq \sigma)$, such that $M_0 = 0$, $M_\alpha = \bigcup_{\beta < \alpha} M_\beta$ for each limit ordinal $\alpha \leq \sigma$, $M_\sigma = M$, and for each $\alpha < \sigma$, $M_{\alpha + 1}/M_\alpha \cong C_\alpha$ for some $C_\alpha \in \mathcal C$. The ordinal $\sigma$ is the \emph{length} of the $\mathcal C$-filtration $\mathcal M$. 

For example, if $\kappa \geq \aleph_0$ and $\mathcal C _\kappa$ denotes the class of all $\leq \kappa$-presented flat Mittag-Leffler modules, then $\Filt {\mathcal C _\kappa} = \mathcal F \mathcal M _\kappa$ is the class of all $\kappa$-restricted flat Mittag-Leffler modules.

\bigskip
\subsection{Direct limits}\label{dirl} Let $(I,\leq)$ be an (upper) directed poset. A covariant functor $\mathfrak C$ from the category $(I,\leq)$ to $\rmod R$ is called an \emph{$I$-direct system of modules}. Equivalently, $\mathfrak C$ may be viewed as a diagram, $\mathcal C$, in the category $\rmod R$, whose objects are modules $C_i$ ($i \in I$), and morphisms are $f_{ji} \in \Hom R{C_i}{C_j}$ for $i \leq j \in I$ which satisfy the identities $f_{ii} = \hbox{id}_{C_i}$ and $f_{kj}f_{ji} = f_{ki}$ for all $i \leq j \leq k \in I$. 

The colimit of the diagram $\mathcal C$ in the category $\rmod R$ is called the \emph{direct limit} of $\mathcal C$ and denoted by $\varinjlim \mathcal C$. In more detail, the colimit is a \emph{cocone} $(C, f_i (i \in I))$ (= a module $C \in \rmod R$ together with morphisms $f_i \in \Hom R{C_i}C$ such that $f_i = f_j f_{ji}$ for all $i \leq j \in I$) possessing the following universal property: for each cocone $(C^\prime, f^\prime_i (i \in I))$ there is a unique homomorphism $\varphi: C \to C^\prime$ such that $\varphi f_i = f^\prime _i$ for each $i \in I$. We will also use the notation $C = \varinjlim_{i \in I} C_i$.

There is a useful presentation of the direct limit as a factor of the direct sum $\bigoplus_{i \in I} C_i$,  
$$0 \to K \overset{\mu}\hookrightarrow \bigoplus_{i \in I} C_i \overset{\pi}\to C = \varinjlim C_i \to 0,$$ 
where $K$ is the submodule of $\bigoplus_{i \in I} C_i$ generated by $\{ x - f_{ji}(x) \mid x \in C_i, i \leq j \in I \}$. Moreover, $\pi \restriction C_i = f_i$ for each $i \in I$, and $\mu$ is a pure embedding. In fact $\mu$ is even locally split, see e.g.\ \cite[4.3]{PPT}. 

There is a convenient way of checking that a cocone $(C^\prime, f^\prime_i (i \in I))$ is a direct limit of the diagram $\mathcal C$. This is the case, iff the following two {\lq\lq}internal{\rq\rq} conditions are satisfied: 

\begin{itemize}
\item[(C1)] $\bigcup_{i \in I} \im f^\prime_i = C^\prime$, and 
\item[(C2)] $\Ker {f^\prime_i} = \bigcup_{i \leq j \in I} \Ker {f_{ji}}$. 
\end{itemize}

Indeed, both conditions hold for the direct limit $(C, f_i (i \in I))$, and condition (C1) is equivalent to the surjectivity of the homomorphism $\varphi$, while (C2) to its injectivity (cf.\ \cite[2.3]{GT}). Notice that condition (C2) implies that if all the morphisms $f_{ji}$ ($i \leq j \in I$) in $\mathcal C$ are monic, then so are all the $f_i$ ($i \in I$). Also, if $I^\prime$ is a $\leq$-cofinal subset of $(I,\leq)$, then conditions (C1) and (C2) hold for the {\lq\lq}restricted{\rq\rq} cocone $(C^\prime, f^\prime_i (i \in I^\prime))$, whence also $C^\prime \cong \varinjlim_{i \in I^\prime} C_i$.

Of particular interest is the case when $I = \sigma$ is a limit ordinal, $\leq$ is the ordinal ordering on $I$, $f_{i+1,i}$ is an inclusion for each $i < \sigma$, and $C_i = \bigcup_{j < i} C_j$ for each limit ordinal $i < \sigma$. Then $( C_i \mid i < \sigma)$ is a \emph{continuous chain} of modules. In this case $C = \bigcup_{i < \sigma} C_i$.        

One can proceed and define \emph{homomorphisms between $I$-direct systems} of modules $\mathcal C$ and $\mathcal C ^\prime$ as systems of morphisms $( \phi_i \mid i \in I )$ such that $\phi_i \in \Hom R{C_i}{C^\prime_i}$ and $\phi_j f_{ji} = f^\prime_{ji} \phi_i$ for all $i \leq j \leq k \in I$. Then $\varinjlim \phi_i : \varinjlim \mathcal C \to \varinjlim \mathcal C^\prime$ is defined by 
$\varinjlim \phi_i (c) = (f^\prime_k \phi_k (c_k))$ for each $c \in C$ such that $c = f_k(c_k)$ for some $k \in I$ and $c_k \in C_k$. 

In this way, $\varinjlim$ defines a functor from the category of $I$-inverse systems of modules to $\rmod R$. This functor is well-known to be exact, see e.g.\ \cite[5.33]{R}.

\bigskip
\subsection{Inverse limits}\label{invl} Let $(I,\leq)$ be a directed poset. A contravariant functor $\mathfrak D$ from the category $(I,\leq)$ to $\rmod R$ is called an \emph{$I$-inverse system} of modules. Equivalently, $\mathfrak D$ may be viewed as a diagram, $\mathcal D$, in the category $\rmod R$, whose objects are modules $D_i$ ($i \in I$), and morphisms are $f_{ij} \in \Hom R{C_j}{C_i}$ for $i \leq j \in I$ that satisfy the identities $f_{ii} = \hbox{id}_{C_i}$ and $f_{ij}f_{jk} = f_{ik}$ for all $i \leq j \leq k \in I$. 

The limit of the diagram $\mathcal D$ in the category $\rmod R$ is called the \emph{inverse limit} of $\mathcal D$ and denoted by $\varprojlim \mathcal D$. The limit is a \emph{cone} $(D, g_i (i \in I))$ (= a module $D \in \rmod R$ and morphisms $g_i \in \Hom RD{D_i}$ such that $g_i = g_{ij}g_j$ for all $i \leq j \in I$) with the following universal property: for each cone $(D^\prime, g^\prime_i (i \in I))$ there is a unique homomorphism $\psi: D^\prime \to D$ such that $g_i \psi = g^\prime _i$ for each $i \in I$. We will also use the notation $D = \varprojlim_{i \in I} D_i$.

The module $D$ is (isomorphic to) a particular submodule of the direct product $\prod_{i \in I} D_i$
$$(\ast) \quad D = \{ (d_i)_{i\in I} \in \prod_{i \in I} D_i \mid d_i = g_{ij}(d_j) \hbox{ for all } i \leq j \in I \},$$ 
and for each $i \in I$, $g_i = \pi_i \restriction D$, where $\pi_i$ is the canonical projection of $\prod_{i \in I} D_i$ on to $D_i$.  

From this presentation of the inverse limit, it follows that if $I^\prime$ is a $\leq$-cofinal subset of $(I,\leq)$, then also $D \cong \varprojlim_{i \in I^\prime} D_i$ for the {\lq\lq}restricted{\rq\rq} cone $(D, g_i (i \in I^\prime))$. In particular, if $I$ a countably infinite set, we can w.l.o.g.\ assume that $I = \omega$. An $I$-inverse system of modules with $I = \omega$ is called a \emph{tower} of modules. 

An inverse system $\mathcal D$ is called a \emph{generalized tower} if $I = \sigma$ is a limit ordinal with the ordinal ordering $\leq$, and $\mathcal D = (D_\alpha, g_{\alpha \beta} \mid \alpha \leq \beta < \sigma )$ is a continuous inverse system of modules, that is, $D_\alpha = \varprojlim_{\beta < \alpha} \mathcal D_\beta$ for each limit ordinal $\alpha < \sigma$. Further, $\mathcal D$ is called a generalized tower of epimorphisms in case all the maps $g_{\alpha \beta}$ ($\alpha \leq \beta < \sigma$) are surjective, or equivalently, $g_{\alpha \alpha + 1}$ is surjective for each $\alpha < \sigma$.

Next, we define \emph{homomorphisms between $I$-inverse systems} of modules $\mathcal D$ and $\mathcal D ^\prime$ as systems of morphisms $( \varphi_i \mid i \in I )$ such that $\varphi_i \in \Hom R{D_i}{D^\prime_i}$ and $\varphi_j f_{ij} = f^\prime_{ij} \varphi_i$ for all $i \leq j \leq k \in I$. Then $\varprojlim \varphi_i : \varprojlim \mathcal D \to \varprojlim \mathcal D^\prime$ is defined by 
$(\varprojlim \varphi_i)((d_i)_{i \in I}) = (\varphi_i (d_i))$ for each $(d_i)_{i \in I} \in D$. 

Thus $\varprojlim$ defines a functor from the category of $I$-inverse systems of modules to $\rmod R$. Since $\varprojlim$ has a left adjoint (provided by the {\lq constant \rq} inverse system functor from $\rmod R$ to the category of all $I$-directed posets), the functor $\varprojlim$ is left exact, see e.g.\ \cite[5.52]{R}. It is not necessarily right exact in $\rmod R$ in general. Here is a simple example demonstrating that right exactness may fail even for towers of short exact sequences: 

\begin{example}\label{simple} Let $p$ be a prime integer. Consider the tower $\mathcal E _p$ of short exact sequences of abelian groups $0 \to \mathbb Z p^i \subseteq \mathbb Z \to \mathbb Z _{p^i} \to 0$ (= free resolutions of the cyclic groups $\mathbb Z _{p^i}$) for $0 < i < \omega$. The connecting triples of morphisms are $(\nu_{i,i+1},\hbox{id}_{\mathbb Z},\pi_{i,i+1})$ ($0 < i < \omega$), where $\nu_{i,i+1} : \mathbb Z p^{i+1} \subseteq \mathbb Z p^{i}$ is the inclusion and $\pi_{i,i+1} : \mathbb Z _{p^{i+1}} \to \mathbb Z _{p^i}$ the projection modulo the socle of $\mathbb Z _{p^{i+1}}$. Then $\varprojlim \mathcal E _p$ is the sequence $0 \to 0 \to \mathbb Z \to \mathbb J _p \to 0$ where $\mathbb J _p = \varprojlim \mathbb Z _{p^{i}}$ is the (uncountable) group of all $p$-adic integers. So $\varprojlim \mathcal E _p$ is not right exact. 
\end{example}    

\medskip
If a cone $(D, g_i(i \in I))$ is the inverse limit of $\mathcal D$, then the following {\lq\lq}internal{\rq\rq} condition (D1) -- dual to the condition (C1) above -- holds

\begin{itemize}
\item[(D1)] $\bigcap_{i \in I} \Ker {g_i} = 0$. 
\end{itemize}

Indeed, in the notation above, for a cone $(D^\prime, g^\prime_i(i \in I))$, condition (D1) is equivalent to the injectivity of the homomorphism $\psi$. However, no {\lq\lq}internal{\rq\rq} condition is known to be equivalent to the surjectivity of $\psi$ in general. Of course, we can formally dualize condition (C2) as

\begin{itemize}
\item[(D2)] $\im g_i = \bigcap_{i \leq j \in I} \im g_{ij}$. 
\end{itemize} 

Notice that (D2) implies that if all the morphisms $g_{ij}$ ($i \leq j \in I$) in $\mathcal D$ are surjective, then so are all the $g_i$ ($i \in I$).   

\medskip
We will briefly discuss condition (D2) for countable inverse systems of modules. By the above, we can w.l.o.g.\ assume that $I = \omega$, that is, $\mathcal D$ is a tower. For each $i < \omega$, let $B_i = \bigcap_{i \leq j < \omega} \im g_{ij}$ and $h_i = g_{i,i+1} \restriction B_{i+1}$. Clearly, $(D2)$ is equivalent to the surjectivity of all the $h_i$ ($i < \omega$).  

Let us restrict further to the particular case of a tower formed by an iteration of a single endomorphism. That is, we consider $A \in \rmod R$ and $f \in \End _R(A)$, and let $D_i = A$ and $g_{i,i+1} = f$ for all $i < \omega$. Then $(D2)$ holds, iff $f(B) = B$, where $B = \bigcap_{i < \omega} f^i(A)$. The latter clearly holds when $f$ is surjective, and it is easy to see that is also holds when $f$ is injective. However, it may fail for a general endomorphism $f$. To demonstrate this fact, we recall the classic construction of totally projective modules due to Walker (cf.\ \cite[1.7]{BS}):

\begin{example}[Walker's towers]\label{Walker} Let $R$ be a discrete valuation domain with a prime element $p \in R$, and $\beta$ be any infinite ordinal. The module $A$ is defined by generators and relations as follows: the generators are labeled by finite sequences $\beta\beta_1\dots\beta_n$ of ordinals such that $\beta > \beta_1 > \dots > \beta_n$. The relations are $p.\beta\beta_1\dots\beta_n\beta_{n+1} = \beta\beta_1\dots\beta_n$ and $p . \beta = 0$. The endomorphism $f$ is the multiplication by $p$ on $A$. For each ordinal $\sigma$, we define a submodule $p^\sigma A$ of $A$ by induction: $p^0 A = A, p^{\sigma + 1}A = p (p^{\sigma} A)$, and $p^{\sigma} A = \bigcap_{\rho < \sigma} p^{\rho} A$ for $\sigma$ a limit ordinal.  

For each ordinal $\alpha \leq \beta$, let $S_\alpha$ be the submodule of $A$ generated by (the cosets of) the generators labeled by the sequences $\beta\beta_1\dots\beta_n\alpha$. Then $f^i(A) = p^i A = S_i$ for each $i < \omega$, $B = \bigcap_{i < \omega} f^i(A) = p^{\omega}A = S_\omega$, and $f(B) = p^{\omega+1} A = S_{\omega + 1} \subsetneq B$. Thus condition (D2) fails for the tower formed by an iteration of the endomorphism $f$. 

In fact, for each $\alpha \leq \beta$, $p^\alpha A = S_{\alpha}$, $p^\beta A \cong R/pR$, and $p^{\beta +1} A = 0$, whence $\varprojlim \mathcal D = 0$.  
\end{example}                 

Notice that in the setting of Example \ref{simple}, condition (D2) holds for all the three towers of modules forming the tower of short exact sequences $\mathcal E _p$, but $\mathcal E _p$ is not right exact.  

Clearly, (D2) holds for all towers of epimorphisms, and more in general, for all generalized towers of epimorphisms: then also all the morphisms $g_i$ ($i \in I$) are surjective. However, the latter (and hence (D2)) may fail for uncountable (non-continuous) well-ordered inverse systems of epimorphisms. 

Our example exhibiting the failure is based on the construction of an \emph{Aronszajn tree}, that is, of a tree $T$ of height $\aleph_1$ with no branches of length $\aleph_1$, such that for each $x \in T$, the set of all successors of $x$ in $T$ has cardinality $\aleph_1$, and all levels $T_\alpha$ ($\alpha < \aleph_1$) of $T$ are countable. We refer to \cite[Appendix on Set Theory]{FS} for a construction of such tree.

\begin{example}[Aronszajn's well-ordered inverse systems]\label{Aronszajn}  Let $T$ be an Aronszajn tree. For each $\alpha < \aleph_1$, let $B_\alpha$ be the set of all branches in $T$ of length $\alpha$. For $\alpha \leq \beta < \aleph_1$, we define a map $z_{\alpha \beta } : B_\beta \to B_\alpha$ as the restriction map. That is, $z_{\alpha \beta }$ restricts each branch $\eta \in B_\beta$ to its initial segment in $B_\alpha$. Since all levels of $T$ are countable, and for each $x \in T$, the set of all successors of $x$ in $T$ has cardinality $\aleph_1$, the maps $z_{\alpha \beta }$ are surjective for all $\alpha \leq \beta < \aleph_1$.  

Let $M$ be any module. Let $I = \aleph_1$ with the ordinal ordering $\leq$. For each $\alpha \in I$, let $D_\alpha = M^{(B_\alpha)}$. For $\alpha \leq \beta \in I$, we define an epimorphism $g_{\alpha \beta}: D_\beta \to D_\alpha$ by $g_{\alpha \beta}((y_\eta)_{\eta \in B_\beta}) = (x_\tau)_{\tau \in B_\alpha}$, where for each $\tau \in B_\alpha$, $x_\tau = \sum_{\eta \in B_\beta, z_{\alpha \beta}(\eta) = \tau} y_\eta \in M$. As $T$ has no branch of length $\aleph_1$, ($\ast$) yields that $D = \varprojlim D_\alpha = 0$.       

Moreover, all the $g_{\alpha \beta}$ ($\alpha \leq \beta \in I$) are epimorphisms, while all the $g_\alpha : D \to D_\alpha$ ($\alpha \in I$) are zero morphisms.
\end{example} 

There is, however, an important class of inverse limits of modules where condition (D2) does hold, namely the class of dual inverse systems: 

\begin{example}[Dual inverse systems]\label{dual} Let $R$ be a ring. Denote by $(-)^* = \Hom {\mathbb Z}{-}{\mathbb Q/\mathbb Z}$ the character module duality from $\rmod R$ to $\lmod R$. Notice that for each homomorphism $f$ in $\rmod R$, there is a canonical isomorphism of left $R$-modules $\phi : (\im f)^* \cong \im f^*$.   

For a directed set $(I,\leq)$ and a covariant functor $\mathfrak C : I \to \rmod R$, we define a contravariant functor $\mathfrak D = (-)^* \circ \mathfrak C : I \to \lmod R$. In other words, if $\mathcal C = ( C_i, f_{ji} \mid i \leq j \in I )$ is the $I$-direct system in $\rmod R$ corresponding to $C$, then $\mathcal D = ( C_i^*, f_{ji}^* \mid i \leq j \in I )$ is the $I$-inverse system of left $R$-modules corresponding to $D$. $\mathcal D$ is called the \emph{dual inverse system} of $\mathcal C$.  

Let $(C, f_i (i \in I))$ be the direct limit of $\mathcal C$ in $\rmod R$, so in the notation of \ref{dirl}, we have the short exact sequence $0 \to K \overset{\mu}\hookrightarrow \bigoplus_{i \in I} C_i \overset{\pi}\to C \to 0$ where $\pi \restriction C_i = f_i$ ($i \in I$), and $K$ is generated by the elements of the form $x - f_ji(x)$ where $x \in C_i$ and $i \leq j \in I$. 

Then $0 \to C^* \overset{\pi^*}\to \prod_{i \in I} C_i^* \overset{\mu^*}\to K^* \to 0$ is exact in $\lmod R$, and $g \in \pi^*(C^*)$, iff $g = (g_i)_{I \in I}$ where $g_i \in C_i^*$ and $g_i - g_jf_{ji} = 0$ for all $i \leq j \in I$. The latter equality just says that $g_i = f_{ji}^*(g_j)$. Denoting by $\pi_i$ the restriction to $\pi^*(C^*)$ of the $i$th canonical projection of $\prod_{i \in I} C_i^*$ on to $C_i^*$, we infer that $(\pi^*(C^*), \pi_i (i \in I))$ is the inverse limit of $\mathcal D$ in $\lmod R$. As $\pi \restriction C_i = f_i$ and $\pi^* (x) = (f_i^*(x))_{i \in I}$ for each $x \in C^*$, the inverse limit is isomorphic to the cone $( C^*, f_i^* (i \in I))$.         

Let $i \in I$. Consider the direct system $\mathcal E$ of short exact sequences $0 \to \Ker {f_{ji}} \hookrightarrow C_i \overset{f_{ji}}\to \im {f_{ji}} \to 0$ ($i \leq j \in I$) with the connecting homomorphisms $(\nu_{kj},\hbox{id}_{C_i},\pi_{kj})$ where $\nu_{kj} : \Ker {f_{ji}} \subseteq \Ker {f_{ki}}$ is the inclusion and $\pi_{kj} : \im {f_{ji}} \to \im {f_{ki}}$ the canonical epimorphism, for all $i \leq j \leq k \in I$. By condition (C2) for the direct system $\mathcal C$, $\Ker {f_i}$ is the directed union of its submodules $\Ker {f_{ji}}$ ($i \leq j \leq k \in I$). It follows that $\varinjlim \mathcal E$ is the short exact sequence $0 \to \Ker {f_{i}} \subseteq C_i \overset{f_{i}}\to \im {f_{i}} \to 0$. 

Applying the duality $(-)^*$ to $\mathcal E$ and the isomorphism $\phi$ above, we obtain the inverse system $\mathcal E ^*$ of short exact sequences  $0 \to \im {f_{ji}^*} \hookrightarrow C_i^* \to (\Ker f_{ji})^* \to 0$ ($i \leq j \in I$) with the connecting homomorphisms $(\mu_{jk},\hbox{id}_{C_i^*},\nu_{kj}^*)$ where $\mu_{jk}: \im f_{ji}^* \subseteq \im f_{ki}^*$ is the inclusion. Applying $(-)^*$ to $\varinjlim \mathcal E$, we infer that $\varprojlim \mathcal E ^* \cong (\varinjlim \mathcal E)^*$ is the short exact sequence $0 \to \im {f_i^*}  \hookrightarrow C_i^* \to (\Ker {f_i})^* \to 0$. Thus $\im {f_i^*} = \bigcap_{i \leq j \in I} \im {f_{ji}^*}$, and condition (D2) holds for the cone $(C^*, f_i^* (i \in I))$.   
\end{example} 

\begin{remark}\label{other} (1) The tower of abelian groups $\mathcal D : \dots \to \mathbb Z _{p^{i+1}} \to \mathbb Z _{p^{i}} \to \dots \to \mathbb Z _{p} \to 0$ is a dual inverse system. It is obtained by applying the character module duality to the direct system $\mathcal C : 0 \to \mathbb Z _{p} \subseteq \dots \subseteq \mathbb Z _{p^{i}} \subseteq \mathbb Z _{p^{i+1}} \subseteq \dots$. Here, $p$ is a prime integer, $\mathbb Z _{p^\infty} = \varinjlim \mathbb Z _{p^{i}}$ is the Pr\"{u}fer $p$-group, while $\mathbb J _p \cong (\mathbb Z _{p^\infty})^* = \varprojlim \mathbb Z _{p^{i}}$ the group of all $p$-adic integers. 

More in general, if $\mathcal C$ is any continuous chain of modules, then the dual inverse system $\mathcal D$ is a generalized tower of modules.    

(2) Condition (D2) holds also for other types of dualities: for example, if $(-)^* = \Hom {R}{-}{N}$ where $N$ is a pure-injective module, then for each covariant functor $\mathfrak C : I \to \rmod R$, the functor $\mathfrak D = (-)^* \circ \mathfrak C$ defines an $I$-inverse system of abelian groups that satisfies condition (D2), see \cite[1.7]{H}.    
\end{remark} 

\subsection{Mittag-Leffler conditions}\label{MitL}   

Mittag-Leffler conditions are stabilization conditions for the decreasing chains of images of the inverse system maps: 

\begin{definition}\label{mldef} Let $\mathcal D = (D_i, g_{ij} \mid i \leq j \in I )$ be an inverse system of modules and let $D = \varprojlim \mathcal D = ( D_i, g_i (i \in I) )$ be its inverse limit. 

\begin{enumerate}
\item[(1)] $\mathcal D$ is \emph{Mittag-Leffler}, provided that for each $i \in I$ there exists $i \leq j \in I$, such that $\im {g_{ij}} = \im {g_{ik}}$ for each $j \leq k \in I$. 
\item[(2)] $\mathcal D$ is \emph{strict Mittag-Leffler}, provided that for each $i \in I$ there exists $i \leq j \in I$, such that $\im {g_{ij}} = \im {g_i}$. 
\end{enumerate}
\end{definition}

Since $\im {g_i} \subseteq \im {g_{ik}}$ for each $i \leq k \in I$, each strict Mittag-Leffler inverse system is Mittag-Leffler. For example, if all the $g_{ij}$ ($ i \leq j \in I$) are epimorphisms, then $\mathcal D$ is Mittag-Leffler.

\begin{remark}\label{wellknown} It is easy to see that the two notions coincide for towers of modules: if $\mathcal D$ is a Mittag-Leffler tower, we can w.l.o.g.\ assume that in \ref{mldef}(1), $j = i + 1$, and then for each $d_i \in \im {g_{i,i+1}}$, by induction on $k > i$, find a $d_k \in \im {g_{i,i+1}}$ such that $g_{k-1,k}(d_k) = d_{k-1}$. Thus $d_i \in \im {g_i}$, and $\mathcal D$ is strict Mittag-Leffler. 

However, Example \ref{Aronszajn} presents a well-ordered inverse system $\mathcal D$ of epimorphisms -- hence a Mittag-Leffler inverse system -- whose inverse limit is $0$. So $\mathcal D$ is not strict Mittag-Leffler.
\end{remark}

Let us record another case of coincidence of the two notions from \cite{H}:

\begin{lemma}\label{d2} Assume that the inverse system $\mathcal D$ satisfies condition (D2). Then $\mathcal D$ is Mittag-Leffler, iff it is strict Mittag-Leffler.

In particular, the equivalence holds for all dual inverse systems of modules. 
\end{lemma}
\begin{proof} This is a simple consequence of the set $I$ being (upper) directed: the equalities $\im {g_{ij}} = \im {g_{ik}}$ for each $j \leq k \in I$ imply that $\im {g_{ij}} = \bigcap_{i \leq k \in I} \im g_{ik}$. By condition (D2), the latter intersection equals $\im {g_i}$. 

The final claim follows from Example \ref{dual} (see also Remark \ref{other}(2)).
\end{proof}

The Mittag-Leffler conditions are sufficient to guarantee exactness of the functor $\varprojlim$ at towers of modules. More precisely, let 
$$( \dagger ) \quad 0 \to \mathcal A \to \mathcal B \to \mathcal C \to 0$$ 
be a short exact sequence of generalized towers of modules indexed by a limit ordinal $\sigma$, and let 
$$( \ddagger ) \quad 0 \to \varprojlim \mathcal A \to \varprojlim \mathcal B \to \varprojlim \mathcal C \to 0$$ 
be the left exact sequence obtained by applying the functor $\varprojlim$ to $( \dagger )$.

\begin{lemma}\label{towers}
\begin{itemize}
\item[(i)] Assume that $\sigma$ has cofinality $\omega$ (e.g., ($\dagger$) is a short exact sequence of towers of modules). Then ( $\ddagger$ ) is exact provided that $\mathcal A$ is Mittag-Leffler.
\item[(ii)] Assume that $\mathcal A$ is a generalized tower of epimorphisms. Then $( \ddagger )$ is exact.   
\end{itemize}
\end{lemma}
\begin{proof} (i) W.l.o.g., we can assume that $\sigma = \omega$. Then for a tower $\mathcal D$, $\varprojlim \mathcal D = \Ker f _{\mathcal D}$ where $f_{\mathcal D} : \prod_{i < \omega} D_i \to \prod_{i < \omega} D_i$ is defined by $f ((d_i)_{i < \omega}) = (d_i - g_{i,i+1}(d_{i+1}))_{i < \omega}$. By the Snake Lemma, $( \ddagger )$ is exact if coker $f_{\mathcal A} = 0$. But the latter is known to hold when $\mathcal A$ is Mittag-Leffler (see e.g.\ \cite[3.6]{GT}). 

(ii) Let $\mathcal A = (A_i, f_{ij} \mid i \leq j < \sigma )$, $\mathcal B = (B_i, g_{ij} \mid i \leq j < \sigma )$, and $\mathcal C = (C_i, h_{ij} \mid i \leq j < \sigma )$ be the generalized towers, and $(\nu_i,\pi_i)$ ($i < \sigma$) the maps such that the short exact sequences $0 \to A_i \overset{\mu_i}\to \mathcal B \overset{\pi_i}\to \mathcal C \to 0$ form a generalized tower of short exact sequences. 

We have to prove that $\pi = \varprojlim \pi_i : \varprojlim B_i \to \varprojlim C_i$ is surjective. Let $c = ( c_i \mid i < \sigma ) \in \varprojlim C_i$. By induction on $i < \sigma$, we will define a sequences $b = ( b_i \mid i < \sigma ) \in \prod_{i \in I} B_i$ such that $\pi_i (b_i) = c_i$ for all $i < \sigma$ and $g_{ki}(b_i) = b_k$ for all $k < i < \sigma$. Then $b \in \varprojlim B_i$ and $\pi(b) = c$. 

First, since $\pi_0$ is surjective, there exists $b_0 \in B_0$ such that $\pi_0(b_0) = c_0$. Is $b$ is defined up to an $i < \sigma$, then we define $b_{i+1}$ as follows: we take any $b^\prime_{i+1} \in B_{i+1}$ such that $\pi_{i+1}(b^\prime_{i+1}) = c_{i+1}$. Since $\pi_i g_{i,i+1}(b^\prime_{i+1}) = h_{i,i+1} \pi_{i+1}(b^\prime_{i+1}) = h_{i,i+1}(c_{i+1}) = c_i = \pi_i (b_i)$, we have $b_i - g_{i,i+1}(b^\prime_{i+1}) = \nu_i(a_i)$ for some $a_i \in A_i$. 

Since $\mathcal A$ is a generalized tower of epimorphisms, $f_{i,i+1}$ is surjective. So there exists $a_{i+1} \in A_{i+1}$ such that $\nu_i(a_i) =  \nu_i f_{i,i+1} (a_{i+1}) = g_{i,i+1} \nu_{i+1} (a_{i+1})$. Let $b_{i+1} = b^\prime_{i+1} + \nu_{i+1} (a_{i+1})$. Then $\pi_{i+1}(b_{i+1}) = c_{i+1}$, and $g_{i,i+1}(b_{i+1}) = g_{i,i+1} (b^\prime_{i+1}) + \nu_i(a_i) = b_i$.

If $i < \sigma$ is a limit ordinal, then since the generalized towers $\mathcal B$ and $\mathcal C$ are continuous, letting $b_i = ( b_j \mid j < i )$, we conclude that $\pi_i (b_i) = ( c_j \mid j < i ) = c_i$, q.e.d.                   
\end{proof}

\subsection{Relative Mittag-Leffler and tilting modules}\label{relMLs}

Mittag-Leffler conditions are closely related to Mittag-Leffler modules and their relative versions. In order to make this clear, we require further notions and results from \cite{AH} and \cite{H} that generalize the classic (absolute) case studied in \cite{RG}.

\begin{definition}\label{defML} Let $M$ be a module.  
\begin{enumerate}
\item[(1)] Let $B$ be an $R$-bimodule. Then $M$ is \emph{$B$-stationary} (\emph{strict $B$-stationary}), provided that $M = \varinjlim C_i$ for some direct system $(C_i, f_{ji} \mid i \leq j \in I)$ consisting of finitely presented modules, such that the inverse system $(\Hom R{C_i}B, \Hom R{f_{ji}}B \mid i \leq j \in I)$ is Mittag-Leffler (strict Mittag-Leffler) in $\rmod {\mathbb Z}$.  

\item[(2)] Let $\mathcal B$ be a class of modules. Then $M$ is \emph{$\mathcal B$-stationary} (\emph{strict $\mathcal B$-stationary}), provided $M$ is \emph{$B$-stationary} (\emph{strict $B$-stationary}) for each $B \in \mathcal B$.
\end{enumerate}
\end{definition}

\begin{remark} (1) The notions from \ref{defML} can equivalently be defined by replacing the existential quantifier with the universal one, that is, by replacing {\lq}for some direct system{\rq} with {\lq}for each direct system{\rq}, see \cite{AH}

(2) If $B$ is a pure-injective module, then each $B$-stationary module is strict $B$-stationary, cf.\ Remark \ref{other}(2) and Lemma \ref{d2}.        
\end{remark}

\begin{definition}\label{definable} A class of modules $\mathcal C$ is \emph{definable} provided that it is closed under direct limits, direct products and pure submodules. For a class of modules $\mathcal Q$, we denote by $\Def {(\mathcal Q)}$ the \emph{definable closure} of $\mathcal Q$, which is the least definable class of modules containing $\mathcal Q$. 

Given a definable class $\mathcal C$ of left (right) $R$-modules, we define its \emph{dual definable class} of right (left) $R$-modules, denoted by $\mathcal C ^{\vee}$, as the least definable class of right (left) $R$-modules containing the character modules $C^* = \Hom {\mathbb Z}{C}{\mathbb Q/\mathbb Z}$ of all modules $C \in \mathcal C$. Then $\mathcal C = (\mathcal C ^{\vee})^{\vee}$ for any definable class of left (right) $R$-modules $\mathcal C$, see e.g.\ \cite[2.5]{R2}. 
\end{definition}

\begin{example}\label{fintype} Let $\mathcal S$ be a class of \emph{FP$_2$-modules} (i.e., the modules $M$ possessing a presentation $M \cong F/G$ where $F$ is finitely generated projective, and $G$ is a finitely presented submodule of $F$). Then the class $\mathcal S^\perp = \{ N \in \rmod R \mid \Ext 1RMN = 0 \hbox{ for all } M \in \mathcal S \}$ is definable in $\rmod R$, it's dual definable class in $\lmod R$ being $\mathcal S ^\intercal = \{ N \in \rmod R \mid \Tor 1RMN = 0 \hbox{ for all } M \in \mathcal S \}$. The definable classes of this form are called \emph{of finite type}, see e.g.\ \cite[3.2]{BT2}. 

For a concrete example, assume that $R$ is a right coherent ring. Then FP$_2$-modules coincide with the finitely presented modules. If $\mathcal S$ denotes the class of all finitely presented modules, then $\mathcal S^\perp$ is the definable class of all \emph{fp-injective} modules, and $\mathcal S ^\intercal$ the dual definable class of all flat left $R$-modules. 
\end{example}

A key relation between relative Mittag-Leffler properties and stationarity was proved in \cite[2.11]{H}: 

\begin{theorem}\label{relate} Let $\mathcal Q$ be a definable class of left $R$-modules and $\mathcal B = \mathcal Q ^{\vee}$ be its dual definable class (of right $R$-modules). Let $M$ be a module.
	Then the following conditions are equivalent:
	\begin{itemize}
\item[(i)] $M$ is $\mathcal Q$-Mittag-Leffler.
\item[(ii)] $M$ is $\{ Q \}$-Mittag-Leffler for each $Q \in \mathcal Q$.
\item[(iii)] $M$ is (strict) $Q^{*}$-stationary for each $Q \in \mathcal Q$.
\item[(iv)] $M$ is $\mathcal B$-stationary.
  \end{itemize} 
\end{theorem}

We will also need the following description of flat $\mathcal Q$-Mittag-Leffler modules from \cite[2.6]{HT}. Recall that a system $\mathcal S$ of submodules of a module $M$ is called \emph{$\aleph_1$-dense} provided that each countable subset of $M$ is contained in an element of $\mathcal S$, and $\mathcal S$ is closed under unions of countable chains.   

\begin{theorem}\label{almfree} Let $\mathcal Q$ be any class of left $R$-modules and $M \in \rmod R$. 

Then $M \in \mathcal D _{\mathcal Q}$, iff $M$ has an $\aleph_1$-dense system consisting of countably generated flat $\mathcal Q$-Mittag-Leffler modules.  
 \end{theorem}

Countably presented absolute Mittag-Leffler modules are \emph{pure-projective}, that is, they are direct summands in direct sums of finitely presented modules. Hence countably presented modules in $\mathcal F \mathcal M$ are projective. There are stronger versions of Theorems \ref{relate} and \ref{almfree} available for the absolute case (we refer to \cite[3.14 and 3.19]{GT} and \cite[2.7(1)]{EGPT} for details):   

\begin{theorem}\label{absolut} 
\begin{enumerate}
\item[(1)] The following conditions are equivalent for a module $M$:
	\begin{itemize}
\item[(i)] $M$ is Mittag-Leffler.
\item[(ii)] $M$ has an $\aleph_1$-dense system consisting of countably generated pure-projective modules.
\item[(iii)] Each finite (or countable) subset of $M$ is contained in a countably generated pure-projective submodule which is pure in $M$.
\item[(iv)] $M$ is $\rmod R$-stationary.
  \end{itemize} 
\item[(2)]  The following conditions are equivalent for a module $M$:
  \begin{itemize}
\item[(i)] $M$ is flat Mittag-Leffler.
\item[(ii)] $M$ has an $\aleph_1$-dense system consisting of countably generated projective modules.
\item[(iii)] Each finite (or countable) subset of $M$ is contained in a countably generated projective submodule which is pure in $M$.
\item[(iv)] $M$ is flat and $R$-stationary.
 \end{itemize} 
\item[(3)] Let $\kappa$ be an infinite cardinal and $M \in \mathcal F \mathcal M$ be $\leq \kappa$-generated. Then $M$ is $\leq \kappa$-presented.
\end{enumerate}
\end{theorem}
   
Theorem \ref{relate} concerns only definable classes of modules. But this is not a serious restriction, because of the following general fact from \cite[2.2]{R1}, see also \cite[2.10]{H}:

\begin{proposition}\label{defenough} Let $\mathcal Q$ be any class of left $R$-modules and $M \in \rmod R$. Then $M$ is $\mathcal Q$-Mittag-Leffler, iff $M$ is $\Def {(\mathcal Q)}$-Mittag-Leffler.  
\end{proposition} 

If $\mathcal Q$ is a definable class of modules, then there is a useful criterion in \cite[\S1]{H} for countably presented flat modules to be $\mathcal Q$-Mittag-Leffler, expressed in terms of vanishing of the Ext-functor:

\begin{lemma}\label{herb} Let $M$ be a countably presented flat module, $\mathcal Q$ a definable class of left $R$-modules, and $\mathcal B = \mathcal Q ^{\vee}$ the dual definable class in $\rmod R$. Then $M$ is $\mathcal Q$-Mittag-Leffler, if and only if $\Ext 1RMB = 0$ for all $B \in \mathcal B$. 
\end{lemma}

\begin{example}\label{fintyp} Let $\mathcal S$ be a class of FP$_2$-modules. By Example \ref{fintype}, we can take $\mathcal Q = \mathcal S ^\intercal$ in Lemma \ref{herb}, whence  $\mathcal B = \mathcal S^\perp$. Thus, if $M$ is a countably presented flat module, then $M \in \mathcal D _{\mathcal Q}$, iff $M \in {}^\perp (\mathcal S ^\perp)$. W.l.o.g., we can assume that $R \in \mathcal S$; then the latter condition is equivalent to $M$ being a direct summand in a module $N$ such that $N$ has a $\mathcal S$-filtration of length $\omega$, see e.g.\ \cite[6.14 and 7.10]{GT}.   
\end{example}

Restricting further the setting of Example \ref{fintyp}, we arrive at the notions of (infinitely generated) tilting modules, and tilting classes:

\begin{definition}\label{tilting} A module $T$ is \emph{tilting}, in case it satisfies the following three properties:
\begin{itemize}
\item[(T1)] $T$ has finite projective dimension, 
\item[(T2)] $\Ext iRT{T^{(I)}} = 0$ for all sets $I$,
\item[(T3)] there exists a finite exact sequence $0 \to R \to T_0 \to T_1 \to \dots \to T_k \to 0$ such that $T_i \in \Add T$ for each $i \leq k$.
\end{itemize}
Here, $\Add T$ denotes the class of all modules isomorphic to direct summands of (possibly infinite) direct sums of copies of the module $T$. 

Let $T$ be a tilting module. Then the class $\mathcal B _T = T^{\perp_\infty}$ is called the \emph{(right) tilting class} induced by $T$, and $\mathcal A _T = {}^\perp \mathcal B _T$ the \emph{left tilting class} induced by $T$. Moreover, $\Add T = \mathcal A _T \cap \mathcal B _T$. 

If $T$ has projective dimension $\leq n$, then $T$ is called \emph{$n$-tilting}, and similarly for the tilting classes $\mathcal A _T$ and $\mathcal B _T$. 

Two tilting modules $T$ and $T^\prime$ are said to be \emph{equivalent} in case $\Add T = \Add {T^\prime}$.
\end{definition}

Tilting classes fit in the setting of classes of finite type due to the following  (see \cite[9.8]{AH} and \cite[13.35 and 13.46]{GT}):

\begin{theorem}\label{tilt} Let $T$ be a tilting module. Let $\mathcal S _T$ be the representative set of all FP$_2$-modules in $\mathcal A _T$. Then $\mathcal B _T = \mathcal S _T ^\perp$, hence $\mathcal B _T$ is a definable class of finite type. 

Let $\mathcal Q _T = \mathcal S _T ^\intercal$ be the dual definable class in $\lmod R$, and $\mathcal C _T$ be the class of all countably presented modules in $\mathcal A _T$. Then $\mathcal A _T = \Filt {\mathcal C _T}$. Moreover, $\mathcal C _T$ coincides with the class of all countably presented $\mathcal Q _T$-Mittag-Leffler modules $M$ such that $M \in \varinjlim \mathcal S _T$.
\end{theorem}

As all extensions of the modules in $\mathcal D _\mathcal Q$ are pure, the classes $\mathcal D _\mathcal Q$ are closed under extensions (and transfinite extensions) for each class $\mathcal Q \subseteq \lmod R$. This is also true of the classes $\mathcal A _T$ from Theorem \ref{tilt}, because $\mathcal A _T = \Filt {\mathcal C _T}$. However, the class of all (not necessarily flat) Mittag-Leffler modules need not be closed under extensions in general: 

\begin{example}\label{nML} Let $R = \mathbb Z$ and $p$ a prime integer. By the Krull-Remak-Schmidt-Azumaya theorem \cite[12.6]{AF}, the pure-projective abelian $p$-groups are isomorphic to direct sums of copies of the groups $\mathbb Z _{p^n}$ ($0 < n < \omega$). In particular, each countable Mittag-Leffler abelian group contains no non-zero element of infinite $p$-height. 

For each $n < \omega$, let $g_n = 1 + \mathbb Z p^{n+1} \in \mathbb Z _{p^{n+1}}$. Consider the short exact sequence 
$$0 \to \mathbb Z _p \overset{\nu}\to \bigoplus_{0 < n < \omega} \mathbb Z _{p^{n+1}}/I \overset{\pi}\to \bigoplus_{0 < n < \omega} \mathbb Z _{p^{n}} \to 0$$
where $I =  \bigoplus_{0 < n < \omega} \mathbb Z (p^n g_n - p^{n+1} g_{n+1})$, $\nu(g_0) = p g_1 + I$, and for each $0 < n < \omega$, $\pi(g_n + I) = g_{n-1}$.

Then the outer terms of this short exact sequence, $\mathbb Z _p$ and $\bigoplus_{0 < n < \omega} \mathbb Z _{p^{n}}$, are countable pure-projective abelian groups. However, the middle term $A = \bigoplus_{0 < n < \omega} \mathbb Z _{p^{n+1}}/I$ is not Mittag-Leffler. Indeed, $0 \neq \nu(g_0) = p (g_1 + I) = p^2 (g_2 + I) = \dots = p^n (g_n + I) = \dots$, so $\nu(g_0)$ is a non-zero element of $A$ of infinite $p$-height. This shows that the class of all (countable) Mittag-Leffler abelian groups is not closed under extensions.   
\end{example}

\medskip
\begin{remark}\label{mlxclosure} In the setting of Theorem \ref{tilt}, we also have $\mathcal A _T \subseteq \varinjlim \mathcal S _T = {}^\intercal (\mathcal S _T ^\intercal) = \varinjlim \mathcal A _T$ (cf.\ \cite[8.40]{GT}), so the class $\varinjlim \mathcal S _T$ is always closed under transfinite extensions. However, the characterization of the countably presented modules from $\mathcal A _T$ as the $\mathcal Q _T$-Mittag-Leffler modules in $\varinjlim \mathcal S _T$ from Example \ref{fintyp} does not extend to arbitrary modules in $\mathcal A _T = \Filt {\mathcal C _T}$. 

For example, if $R$ is any non-right perfect ring and $T = R$, then $\mathcal A _ T = \mathcal P _0$, $\mathcal Q _T = \lmod R$ and $\varinjlim \mathcal S _T = \mathcal F _0$, but the class of all $\mathcal Q _T$-Mittag-Leffler modules in $\varinjlim \mathcal S _T$ is just the class $\mathcal F \mathcal M (\supsetneq \mathcal P _0)$. Notice that this example also shows that the criterion for countably presented flat modules to be $\mathcal Q$-Mittag-Leffler from Lemma \ref{herb} does not extend to all flat modules (here, $\mathcal D _{\mathcal Q} = \mathcal F \mathcal M$, while $^\perp \mathcal B = \mathcal P _0$).               
\end{remark}

\section{Approximations}\label{approx}

Precovering classes $\mathcal C$ of modules are important for extending classical homological algebra to more refined settings. Classically, one uses projective resolutions of modules. In the refined setting, one deals with $\mathcal C$-resolutions obtained by iterations of $\mathcal C$-precovers. This is one of the themes of relative homological algebra \cite{EJ}.  

\begin{definition} Let $\mathcal C$ be a class of modules and $M \in \rmod R$. A homomorphism $f : C \to M$ with $C \in \mathcal C$ is called a \emph{$\mathcal C$-precover} of $M$ (or a \emph{right $\mathcal C$-approximation} of $M$) provided that for each homomorphism $f^\prime : C^\prime \to M$ with $C^\prime \in \mathcal C$, there exists a homomorphism $g : C^\prime \to C$ such that $f^\prime = f g$. The $\mathcal C$-precover $f$ is a \emph{$\mathcal C$-cover}, provided that $f$ is \emph{right minimal}, i.e., if each $g \in \End {_R(C)}$ such that $f = fg$ is an automorphism of $C$.   

$\mathcal C$ is called a \emph{precovering (covering)} class, if each module $M \in \rmod R$ has a $\mathcal C$-precover ($\mathcal C$-cover).
\end{definition}   
  
It is well-known that the class $\mathcal P_0$ is precovering for any ring $R$, and $\mathcal P_0$ is covering, iff $R$ is a right perfect ring. However, the class $\mathcal F _0$ is covering for every ring $R$, \cite{BEE}. Our goal in this section is to investigate precovering properties of the intermediate classes of the restricted, and the relative, flat Mittag-Leffler modules. Our presentation follows \cite{BT1}, \cite{S}, and \cite{ST}. 

\subsection{Approximations by restricted Mittag-Leffler modules}\label{aresMLs}

We start with recalling a remarkable general result from \cite[2.15]{SaSt} (see also \cite[7.21]{GT}):

\begin{theorem}\label{filt} Let $\mathcal S$ be any set of modules. Then the class $\Filt {\mathcal S}$ is precovering.  
\end{theorem} 

Theorem \ref{filt} of course includes the case of restricted flat Mittag-Leffler modules:

\begin{corollary}\label{restr} The class $\mathcal F \mathcal M _\kappa$ is precovering for each infinite cardinal $\kappa$. 
\end{corollary}

The class  $\mathcal F \mathcal M _{\aleph_0} = \mathcal P_0$ is not covering when $R$ is not right perfect, \cite[28.4]{AF}. It is conjectured that the same holds for the classes $\mathcal F \mathcal M _\kappa$ when $\kappa > \aleph_0$. In fact, there is a much more general conjecture due to Enochs:

\begin{conjecture}\label{enochs}[{\bf Enochs' Conjecture}] Let $\mathcal C$ be a precovering class of modules. Then $\mathcal C$ is covering, iff $\mathcal C = \varinjlim \mathcal C$.
\end{conjecture}

Enochs' Conjecture is still open in general, but has been proved in a number of particular cases (e.g., for all left tilting classes in \cite[5.2]{AST}; note that all left tilting classes are precovering by Theorems \ref{tilt} and \ref{filt}). Recently, an important case of the conjecture was proved to be \emph{consistent with ZFC} in \cite{BSa}: it holds in any extension of ZFC satisfying G\"{o}del's Axiom of Constructibility (V = L):

\begin{theorem}\label{VL} Assume V = L. Let $\mathcal S$ be any set of modules. Then the class $\Filt {\mathcal S}$ is covering, iff it is closed under direct limits.  
\end{theorem} 

Since $\mathcal F \mathcal M _\kappa$ is not closed under direct limits for any $\kappa \geq \aleph_0$ in case $R$ is not right perfect, we have

\begin{corollary}\label{consist} Assume that V = L. Let $R$ be a non-right perfect ring. Then the class $\mathcal F \mathcal M _\kappa$ is not covering for any infinite cardinal $\kappa$. 
\end{corollary}

\subsection{Approximations by relative Mittag-Leffler modules}\label{arelMLs}  

\emph{For the rest of this section, we will again assume that $R$ is a non-right perfect ring}. Then the setting of relative Mittag-Leffler modules is quite different from the restricted ones: we will see that for any class of left modules $\mathcal Q$, the class of relative flat Mittag-Leffler modules $\mathcal D _{\mathcal Q}$ is precovering only in the boundary case of $\mathcal D _{\mathcal Q} = \mathcal F _0$. This will follow from the next two lemmas that were originally proved in more general forms in \cite{ST} and \cite{S}. First, we need further notation.

\begin{definition}\label{setting} Let $N$ be a countably presented flat non-projective module. (Such modules exist, because $R$ is not right pefect; they are called \emph{Bass modules}). Since $\mathcal F \mathcal M _{\aleph_0} = \mathcal P _0$, necessarily $N \notin \mathcal F \mathcal M$. As $N$ is a countable direct limit of finitely generated free modules, there is a chain 
$$( \ast \ast ) \quad F_0 \overset{h_0}\to F_1 \to \dots \overset{h_{i-1}}\to F_i \overset{h_i}\to F_{i+1} \overset{h_{i+1}}\to \dots$$
where $F_i$ is a finitely generated free module for each $i < \omega$ such that $N \cong \varinjlim_{i < \omega} F_i$.  

Let $\kappa$ be an infinite cardinal and $T_\kappa$ be the \emph{tree on $\kappa$} consisting of all finite sequences of ordinals $< \kappa$. That is, each $\tau \in T_\kappa$ is a map $\tau : n \to \kappa$ for some $n < \omega$. The partial order on $T_\kappa$ is by restriction, so if $\tau, \rho \in T_\kappa$, then $\tau \leq \rho$, iff  $\rho \restriction n = \tau$. For each $\tau \in T_\kappa$, $\ell(\tau)$ will denote the length of $\tau$. 

Let $B_\kappa$ denote the set of all branches of $T_\kappa$. Notice that $\card {T_\kappa} = \kappa$, and $\card {B_\kappa} = \kappa^\omega$ (because branches in $T_\kappa$ correspond to $\omega$-sequence of ordinals $< \kappa$).
\end{definition} 

The following lemma is a special instance of Lemma \cite[5.6]{ST}:

\begin{lemma}\label{tree} There exists a module $L_\kappa \in \mathcal F \mathcal M$ (called the \emph{tree module} for $\kappa$), such that $L_\kappa$ contains a free submodule $D_\kappa$ of rank $\kappa$, and $L_\kappa/D_\kappa \cong N^{(\kappa^\omega)}$. 
\end{lemma}
\begin{proof} Let $D_\kappa = \bigoplus_{\tau \in T_\kappa} F_{\ell(\tau)}$, and $P_\kappa = \prod_{\tau \in T_\kappa} F_{\ell(\tau)}$. Since $\card {T_\kappa} = \kappa$, $D_\kappa$ is a free module of rank $\kappa$

For each $\nu \in B_\kappa$, $i < \omega$, and $x \in F_i$, we define $x_{\nu i} \in P_\kappa$ by $\pi_{\nu \restriction i} (x_{\nu i}) = x$, $\pi_{\nu \restriction j} (x_{\nu i}) = h_{j-1}\dots h_i(x)$ for all $i < j < \omega$, and $\pi_\tau (x_{\nu i}) = 0$ otherwise. Here $\pi_\tau \in \Hom {R}{P_\kappa}{F_{\ell(\tau)}}$ denotes the $\tau$th projection for each $\tau \in T_\kappa$.

Let $Y_{\nu i} = \{ x_{\nu i} \mid x \in F_i \}$. Then $Y_{\nu i}$ is a submodule of $P_\kappa$ isomorphic to $F_i$ via the assignment $x \mapsto x_{\nu i}$. Let $X_{\nu i} = \sum_{j \leq i} Y_{\nu j}$. Then $X_{\nu i} \subseteq X_{\nu, i+1}$, and $X_{\nu i} = \bigoplus_{j < i} F_j \oplus Y_{\nu i} \cong \bigoplus_{j \leq i} F_j$. Let $X_\nu = \bigcup_{i < \omega} X_{\nu i}$, and $L_\kappa = \sum_{\nu \in B_\kappa} X_\nu$. 

Notice that $X_\nu \cong \bigoplus_{i < \omega} F_i$ for each $\nu \in B_\kappa$. Indeed, the inclusion $X_{\nu i} \subseteq X_{\nu, i+1}$ splits, as the short exact sequence 
$0 \to Y_{\nu i} \overset{p}\hookrightarrow F_i \oplus Y_{\nu, i+1} \overset{q}\to F_{i+1} \to 0$ splits, where $p (x_{\nu i}) = x + (g_i(x))_{\nu, i+1}$, and $q (x + y_{\nu, i+1}) = y - g_i(x)$.  

Further, for each $\nu \in B_\kappa$, $N \cong (X_\nu + D_\kappa)/D_\kappa$, as for each $i < \omega$, we can define $f_i : F_i \to (X_\nu + D_\kappa)/D_\kappa$ by $f_i(x) = x_{\nu i} + D_\kappa$, and $((X_\nu + D_\kappa)/D_\kappa, f_i \mid i \in I)$ is the direct limit of the direct system ($\ast \ast$). 

Moreover, each element of $X_\nu$ is a sequence in $P_\kappa$ whose $\tau$th component is zero for all $\tau \notin \{ \nu \restriction i \mid i < \omega \}$, so the modules $((X_\nu + D_\kappa)/D_\kappa \mid \nu \in B_\kappa) $ are independent. It follows that $L_\kappa/D_\kappa = \bigoplus_{\nu \in  B_\kappa} (X_\nu + D_\kappa)/D_\kappa \cong N^{(B_\kappa)}$.

For each countable subset $C = \{ \nu_i \mid i < \omega \}$ of $B_\kappa$, the module $X_C = \sum_{\nu \in C} X_\nu$ is isomorphic to a countable direct sum of the $F_i$s. Indeed, $X_C = \bigcup_{i < \omega} X_{C_i}$, where $X_{C_i} = \sum_{j \leq i} X_{\nu_j}$ is a direct summand in $X_{C_{i+1}}$, with the complementing direct summand isomorphic to a countable direct sum of the $F_i$s. It follows that the set $\mathcal S$ of all $X_C$, where $C$ runs over all countable subsets of $B_\kappa$, is an $\aleph_1$-dense system of submodules of $L_\kappa$ consisting of countably generated free modules. By Theorem \ref{absolut}(2), $L_\kappa \in \mathcal F \mathcal M$. 
\end{proof}

Next we need a special case of Lemma \cite[3.2]{S}:   

\begin{lemma}\label{noprec} Let $\mathcal Q$ be a class of left $R$-modules. Assume there exists a Bass module $N \notin \mathcal D _{\mathcal Q}$. Then $N$ does not have a $\mathcal D _{\mathcal Q}$-precover.  
\end{lemma}
\begin{proof} Let $\pi : M \to N$ be a $\mathcal D _{\mathcal Q}$-precover of $N$. Since $\mathcal P_0 \subseteq \mathcal D _{\mathcal Q}$, $\pi$ is surjective, and we have a short exact sequence $0 \to K \overset{\mu}\to M \overset{\pi}\to N \to 0$. Let $\kappa$ be an infinite cardinal such that $\kappa^\omega = 2^\kappa$, $\card R \leq \kappa$, and $\card K \leq 2^\kappa$ (e.g., let $\lambda_0 = \card R + \card K + \aleph_0$, $\lambda_{i + 1} = 2^{\lambda_i}$ and $\kappa = \sup_{i < \omega} \lambda_i$).  

By Lemma \ref{tree}, there is a short exact sequence involving the tree module $L_\kappa$ as follows: $0 \to R^{(\kappa)} \hookrightarrow L_\kappa \to N^{(2^\kappa)} \to 0$.

Consider the group homomorphism $\Ext 1R{L_\kappa}{\mu}$. It takes a short exact sequence $0 \to K \to X \to L_\kappa \to 0$ to its pushout with $\mu$. The pushout yields a commutative diagram with exact rows and columns as follows:
$$
\begin{CD}
@.       0@.       0              @.       @.
\\
@.       @VVV      @VVV           @.       @.
\\
0@>>>    K@>{\mu}>> M@>{\pi}>>    N@>>>    0
\\
@.       @VVV      @V{\delta}VV           @|       @.
\\
0@>>>    X@>{\subseteq}>>     Y@>{\rho}>>          N@>>>    0
\\
@.       @V{\sigma}VV      @V{\tau}VV    @.       @.
\\
0@>>>    {L_\kappa}@=     {L_\kappa}   @.       @.
\\
@.       @VVV      @VVV           @.       @.
\\
@.      0@.       0.              @.       @.
\end{CD}
$$
Assume that the middle vertical sequence splits, so there exists $\nu \in \Hom R{L_\kappa}Y$ such that $\tau \nu = 1_{L_\kappa}$. Since $\pi$ is a  
$\mathcal D _{\mathcal Q}$-precover of $N$ and $L_\kappa \in \mathcal F \mathcal M \subseteq \mathcal D _{\mathcal Q}$, there exists $\phi \in \Hom R{L_{\kappa}}M$ such that $\rho \nu = \pi \phi$. Then $\rho \nu = \rho \delta \phi$. Thus $\nu - \delta \phi$ maps into $X$, and 
$\sigma (\nu - \delta \phi) = \tau (\nu - \delta \phi) = \tau \nu - 0 = 1_{L_\kappa}$. Hence also the left vertical sequence splits. 
This proves that the group homomorphism $\Ext 1R{L_\kappa}{\mu}$ is monic.

Consider the commutative diagram with exact rows 

$$\begin{CD}
\Hom R{R^{(\kappa)}}K @>{h}>> {\Ext 1RNK}^{2^\kappa} @>{f}>> \Ext 1R{L_\kappa}K \\
 @V{\Hom R{R^{(\kappa)}}{\mu}}VV  @V{{\Ext 1RN{\mu}}^{2^\kappa}}VV  @V{\Ext 1R{L_\kappa}{\mu}}VV \\
\Hom R{R^{(\kappa)}}M   @>>>   {\Ext 1RNM}^{2^\kappa} @>>> {\Ext 1R{L_\kappa}M}.
\end{CD}$$

 Since $\Ext 1R{L_\kappa}{\mu}$ is monic, $\Ker {{\Ext 1RN{\mu}}^{2^\kappa}} \subseteq \Ker f = \im h$. If $\Ker {\Ext 1RN{\mu}} \neq 0$, then $\Ker {{\Ext 1RN{\mu}}^{2^\kappa}}$ has cardinality $\geq 2^{2^\kappa}$, while $\im h$ has cardinality less than or equal to $\card \Hom R{R^{(\kappa)}}K = \card K^\kappa \leq 2^\kappa$. Thus, also $\Ext 1RN{\mu}$ is monic, which implies that $\Hom RN{\pi}$ is surjective. In particular, $\pi$ splits, so $N \in \mathcal D _{\mathcal Q}$, a contradiction.        
\end{proof} 

Now we can proceed as in \cite[Theorem 2.6]{BT1}:

\begin{theorem}\label{boundary} Let $\mathcal Q$ be any class of left $R$-modules. Then the class $\mathcal D _{\mathcal Q}$ is precovering, iff $\mathcal D _{\mathcal Q}$ is deconstructible, iff $\mathcal D _{\mathcal Q} = \mathcal F _0$.
\end{theorem} 
\begin{proof} The class $\mathcal F _0$ is deconstructible over any ring, and each deconstructible class is precovering by Theorem \ref{filt}.  

Assume there exists $M \in \mathcal F _0 \setminus \mathcal D _{\mathcal Q}$. As $M \in \mathcal F _0$, $M$ is a direct limit of a direct system of finitely generated free modules $( F_i \mid i \in I)$ for a directed set $(I,\leq)$. By \cite[3.11]{GT}, there exists a countable chain $i_0 < i_1 < \dots i_n < i_{n+1} < \dots$ of elements of $I$ and a countable family $( Q_j \mid j < \omega )$ of elements of $\mathcal Q$ such that $N = \varinjlim_{n < \omega} F_{i_n} \notin \mathcal D _{\mathcal Q ^\prime}$ where $\mathcal Q ^\prime = \{ Q_j \mid j < \omega \}$. As $\mathcal D _{\mathcal Q} \subseteq \mathcal D _{\mathcal Q ^\prime}$, $N$ is a Bass module which is not contained in $\mathcal D _{\mathcal Q}$. By Lemma \ref{noprec}, $N$ has no $\mathcal D _{\mathcal Q}$-precover.         
\end{proof}

\begin{example}\label{coherent}  If $\mathcal Q = \{ R \}$, then the modules in $\mathcal D _{\mathcal Q}$ are called \emph{f-projective}, \cite{G}. Recall that a module $M$ is \emph{coherent} in case all finitely generated submodules of $M$ are finitely presented (so a ring $R$ is right coherent, if the regular module $R$ is coherent). 

By \cite[3.5]{BT1}, if $R$ is right coherent, then f-projective modules are exactly the flat coherent modules, whence $\mathcal D _{\{ R \}} = \mathcal F _0$, iff all flat modules are coherent. By \cite[3.6]{BT1}, the latter holds for each coherent domain $R$.      
\end{example}

\section{Zariski locality}\label{Zariski}

\emph{All rings in this section are commutative}. By a classic theorem of Grothendieck, if $R$ is a ring, then the category $\Qcoh X$ of all quasi-coherent sheaves on the affine scheme $X = \Spec R$ is equivalent to the module category $\rmod R$. 

For general schemes $X$, $\Qcoh X$ can be represented as a category consisting of \emph{quasi-coherent $\mathcal R$-modules} $\mathcal M = ( M(U), f_{UV} \mid V \subseteq U \subseteq X, U, V \hbox{ open affine })$ over the structure sheaf of rings $\mathcal R = ( R(U) \mid U \hbox{ open affine subset of } X )$ as follows: for every open affine subset $U \subseteq X$, $M(U)$ is an $R(U)$-module, and for each pair of open affine subsets $V \subseteq U \subseteq X$, $f_{UV}: M(U) \to M(V)$ is an $R(U)$-homomorphism such that 

\begin{itemize}
\item $\mbox{id}_{R(V)} \otimes f_{UV}$ is an $R(V)$-isomorphism, and
\item if $W$ is an open affine subset of $V$, then $f_{UV} f_{VW} = f_{UW}$.
\end{itemize}

One does not need all the open affine subsets of $X$ to represent $\Qcoh X$ in this way: an open affine covering $\mathcal C$ of $X$ is sufficient for this purpose, cf.\ \cite[\S2]{EE}. Here, a set $\mathcal C$ of open affine subsets of $X$ is an \emph{open affine covering} of $X$ in case $\mathcal C$ covers both $X$, and all the sets $U \cap V$ where $U$ and $V$ are open affine subsets of $X$. We can view such $\mathcal C$ as a \emph{coordinate system} on $X$ that enables us to represent the geometric object of interest (a quasi-coherent sheaf on $X$) by an algebraic one (a quasi-coherent $\mathcal R$-module). 

\medskip
For any property $P_R$ of modules over a ring $R$, one can use the representation above to extend $P_R$ to a property $P$ of quasi-coherent sheaves on schemes as follows: a quasi-coherent $\mathcal R$-module $\mathcal M$ is \emph{locally $P$-quasi-coherent} in case for each open affine subset $U$ of $X$, the $R(U)$-module $M(U)$ has property $P_{R(U)}$. We will call $P$ the property of quasi-coherent sheaves \emph{induced} by the property $P_R$ of rings $R$. 

For example, if $P_R$ is the property of being a projective $R$-module, then the locally projective quasi-coherent sheaves on a scheme $X$ are exactly the vector bundles on $X$, \cite{D}.  

\medskip  
Of course, we are interested in those properties $P$ that are independent of the choice of a coordinate system, so that they can be checked using \emph{any} open affine covering $\mathcal C$ of $X$. Then an $\mathcal R$-module $\mathcal M$ is locally $P$-quasi-coherent, whenever the $R(U)$-module $M(U)$ has property $P_{R(U)}$ for each $U \in \mathcal C$. Such properties $P$ of quasi-coherent sheaves on $X$ are called \emph{Zariski local}. We will also say that the notion of a locally $P$-quasi-coherent sheaf is Zariski local, or \emph{affine local}.  

\medskip
Our goal here is to prove Zariski locality for the various notions of locally $P$-quasi-coherent sheaves induced by classes of restricted and flat relative Mittag-Leffler modules, and by tilting modules. The definition of a locally $P$-quasi-coherent sheaf given above will be sufficient to achieve this goal in the case of vector bundles (in section \ref{advb}), their generalizations to $\kappa$-restricted Drinfeld vector bundles (\ref{adres}), and in the case of locally $n$-tilting quasi-coherent sheaves (\ref{zartilt}). 

However, for relative Mittag-Leffler modules (\ref{ftrel}), we will have to impose extra compatibility conditions on the relations among the properties $P_R$ for various rings $R$, and possibly also restrict the type of schemes considered. 

\medskip 
Our main tool for proving the Zariski locality is the following classic lemma \cite[5.3.2]{V} (see also \cite[3.5]{EGT}):

\begin{lemma}\label{acl}[The Affine Communication Lemma] Let $R$ be a ring, $M \in \rmod R$, and $P _R$ be a property of $R$-modules such that
\begin{enumerate} 
\item if $M$ satisfies property $P _R$, then $M[f^{-1}] = M \otimes _R R[f^{-1}]$ satisfies property $P _{R[f^{-1}]}$ for each $f \in R$, and
\item if $R = \sum_{j < m} f_jR$, and the $R[f_j^{-1}]$-modules $M[f_j^{-1}] = M \otimes_R R[f_j^{-1}]$ satisfy property $P _{R[f_j^{-1}]}$ for all $j < m$, then $M$ satisfies property $P _R$.
\end{enumerate}
Then the induced property $P$ of quasi-coherent sheaves on $X$ is Zariski local for every scheme $X$.
\end{lemma}

Notice that for each $f \in R$, the localization in $f$, $\varphi_f : R \to R[f^{-1}]$, is a flat ring homomorphism (that is, $\varphi_f$ makes $R[f^{-1}]$ into a flat $R$-module). Moreover, the ring homomorphism $\varphi_{f_0,...,f_{m-1}} : R \to \prod_{i < m} R[f_i^{-1}]$ is faithfully flat when $R = \sum_{j < m} f_jR$ (that is, $\varphi_{f_0,...,f_{m-1}}$ makes  $\prod_{i < m} R[f_i^{-1}]$ into a faithfully flat $R$-module). So the assumptions of the Affine Communication Lemma are satisfied in case $P$ ascends along flat ring homomorphisms, and descends along faithfully flat ring homomorphisms in the sense of the following definition:

\begin{definition}\label{asdes} Let $\varphi : R \to S$ be a flat ring homomorphism, and $P$ be a property of modules.
\begin{itemize}
\item[(i)] $P$ is said to \emph{ascend} along $\varphi$ if for each $R$-module $M$ with the property $P_R$, the $S$-module $M \otimes_R S$ has the property $P_S$.
\item[(ii)] Assume $\varphi$ is a faithfully flat ring homomorphism. Then $P$ is said to \emph{descend} along $\varphi$ if for each $R$-module $M$, $M$ has the property $P_R$ whenever the $S$-module $M \otimes_R S$ has the property $P_S$. 
\end{itemize}
If $P$ ascends along all flat ring homomorphisms, and descends along all faithfully flat ring homomorphisms, then $P$ is called an \emph{ad-property}.
\end{definition}

In view of Lemma \ref{acl}, in order to prove Zariski locality of a property $P$ of quasi-coherent sheaves on a scheme $X$, it suffices to verify that $P$ is an ad-property. This is the way we will proceed below for the properties arising from Mittag-Leffler conditions. 

\subsection{Vector bundles}\label{advb}

Let us start with the model case of vector bundles going back to \cite[Seconde partie]{RG}. Let $P_R$ be the property of being a projective $R$-module. As mentioned above, a quasi-coherent sheaf on a scheme $X$ is locally $P$-quasi-coherent, iff it is a vector bundle on $X$. 

While the ascent of $P$ is trivial, the descent is a nontrivial fact: First, notice that a module $M$ is projective, iff $M$ is flat Mittag-Leffler and $M$ decomposes into a direct sum of countably presented modules (the if part follows from the fact that countably presented flat Mittag-Leffler modules are projective, cf.\ Theorem \ref{absolut}(2)). 

The tools we have presented so far make it possible to prove the descent of the property of being a flat Mittag-Leffler module (cf.\ \cite[7.33]{GT} or \cite[9.2]{P}):

\begin{lemma}\label{descfML} The property of being a flat Mittag-Leffler module is an ad-property.
\end{lemma}
\begin{proof} As the ascent of projectivity is trivial, the ascent of the property of being a flat Mittag-Leffler module follows by Theorem \ref{absolut}(2). 

To prove the descent, let $\varphi : R \to S$ be a faithfully flat ring homomorphism, and $M \in \rmod R$ be such that $M \otimes _RS \in \rmod S$ is flat and Mittag-Leffler. Viewed as an $R$-module, $M \otimes _RS$ is also flat, because it is isomorphic to $M \otimes _R (S \otimes _S S) \cong (M \otimes _R S) \otimes _S S$ and $S$ is a flat $R$-module (via $\varphi$). So for each short exact sequence $\mathcal E$ of $R$-modules, $\mathcal E \otimes _R (M \otimes _R S)$ is a short exact sequence of $S$-modules. Hence, by faithful flatness of $\varphi$, $\mathcal E \otimes _R M$ is exact in $\rmod R$, whence $M$ is a flat $R$-module. 

Thus $M$ is isomorphic to the direct limit of a direct system $\mathcal C = ( F_i, f_{ji} \mid i \leq j \in I )$ of finitely generated free modules, $M = \varinjlim_{i \in I} F_i$. Applying the functor $\Hom R{-}R$, we obtain the inverse system $\mathcal D = ( \Hom R{F_i}R, \Hom R{f_{ji}}R \mid i \leq j \in I )$. By Theorem \ref{absolut}(2), we have to prove that $M$ is $R$-stationary, i.e., the inverse system $\mathcal D$ is Mittag-Leffler. 

Notice that $M \otimes_R S = \varinjlim_{i \in I} F_i \otimes_R S$, and $M \otimes_R S$ is a Mittag-Leffler $S$-module by assumption. By Theorem \ref{absolut}(1), for each $i \in I$ there exists $j \geq i$, such that $\im {\Hom S{f_{ji}\otimes_R S}S} = \im {\Hom S{f_{ki}\otimes_R S}S}$ for all $j \leq k \in I$. Since $R$ is commutative, there is a natural homomorphism
$\Hom RFR \otimes_R S \to \Hom S{F \otimes_R S}S$; if $F$ is finitely generated and free, then this is even an isomorphism. The faithful flatness of $\varphi$ thus yields that $\im {\Hom R{f_{ji}}R} \otimes_R S = \im {\Hom R{f_{ki}}R} \otimes_R S$ for all $j \leq k \in I$. Again by faithful flatness, we conclude that $\mathcal D$ is a Mittag-Leffler inverse system.
\end{proof}

The fact that projectivity (= the property of being flat Mittag-Leffler, and a direct sum of countably presented modules) descends along faithfully flat ring homomorphisms of commutative rings can now be proved by a technique called \emph{d\'evisage}, \cite[Seconde partie]{RG} (see also \cite[9.6]{P}): As a first step, we deduce from Lemma \ref{descfML} that if $M \otimes_R S$ is a countably generated projective $S$-module, then $M$ is a countably generated projective $R$-module. Then we fix a decomposition of the module $M \otimes_R S \cong \bigoplus_{i \in I} Q_i$ into a direct sum of countably presented projective $S$-modules, and use it to construct by induction on $\alpha$ a continuous chain $( M_\alpha \mid \alpha < \sigma )$ of submodules of $M$ such that $M = \bigcup_{\alpha < \sigma} M_\alpha$, and $M_\alpha \otimes_R S = \bigoplus_{i \in I_\alpha} Q_i$ for a subset $I_\alpha$ of $I$, so that $\card {(I_{\alpha+1} \setminus I_\alpha)} \leq \aleph_0$ for each $\alpha < \sigma$. As $(M_{\alpha +1}/M_\alpha ) \otimes_R S \cong \bigoplus_{i \in I_{\alpha +1} \setminus I_\alpha} Q_i$, the first step yields that $M_{\alpha +1}/M_\alpha$ is a projective module for every $\alpha < \sigma$, whence $M$ is projective. Thus, we obtain 

\begin{theorem}\label{vecbun} The property of being a projective module is an ad-property. Hence the notion of a vector bundle is Zariski local for all schemes.
\end{theorem}  

\subsection{Quasi-coherent sheaves arising from restricted flat Mittag-Leffler modules}\label{adres}

Let $\kappa \geq \aleph_0$ and $P_R$ be the property of being a $\kappa$-restricted flat Mittag-Leffler $R$-module. Let $P$ be the induced property of quasi-coherent sheaves. Then the locally $P$-quasi-coherent sheaves are called \emph{$\kappa$-restricted Drinfeld vector bundles}, cf.\ \cite[p.1423]{EGPT} and \cite[2.1.3]{EGT}.     

Notice that for $\kappa = \aleph_0$, $\mathcal P _0 = \mathcal F \mathcal M _{\aleph_0}$, because countably presented flat Mittag-Leffler modules are just the countably presented projective modules. So by Theorem \ref{vecbun}, the question of whether the property of being a locally $P$-quasi-coherent sheaf is Zariski local has a positive answer in the particular case of $\kappa = \aleph_0$.
 
In \cite[1.1]{EGT}, Theorem \ref{vecbun} was generalized to provide a positive answer for each infinite cardinal $\kappa$. We will now present this result with a simplified proof: 

\begin{theorem}\label{krestric} Let $\kappa$ be an infinite cardinal and $P_R$ be the property of being a $\kappa$-restricted flat Mittag-Leffler $R$-module. Then $P_R$ is an ad-property. Hence the notion of a $\kappa$-restricted Drinfeld vector bundle is Zariski local for all schemes. 
\end{theorem}
\begin{proof} The proof goes along the lines of the proof of Theorem \ref{vecbun} above except for the final part concerning the descent along faithfully flat ring homomorphisms. Here, d\'{e}visage is replaced by a more general technique dealing with filtrations of modules rather than their direct sum decompositions. 

First, we fix a filtration $\mathcal F$ of the module $M \otimes_R S$ by $\leq \kappa$-presented flat Mittag-Leffler $S$-modules witnessing that $M \otimes_R S$ is a $\kappa$-restricted flat Mittag-Leffler $S$-module. The filtration $\mathcal F$ is then enlarged into a family, $\mathcal H$, of $S$-submodules of $M \otimes_R S$ with the following properties: 

\begin{itemize}
\item[(P1)] $\mathcal H$ is a complete distributive sublattice of the complete modular lattice of all $S$-submodules of $M \otimes_R S$, 
\item[(P2)] if $N, N^\prime \in \mathcal H$ satisfy $N \subseteq N^\prime$, then $N^\prime/N$ is a $\kappa$-restricted flat Mittag-Leffler $S$-module,
\item[(P3)] for each $N \in \mathcal H$ and each subset $X$ of $M \otimes_R S$ of cardinality $\leq \kappa$, there exists $N^\prime \in \mathcal H$ such that $N \cup X \subseteq N^\prime$ and $N^\prime/N$ is $\leq \kappa$-presented.
\end{itemize}

Such enlargement of $\mathcal F$ is possible by a general construction known as the Hill lemma, see \cite[7.10]{GT}. Let $\mathcal C _\kappa$ be the class of all $\leq \kappa$-presented flat Mittag-Leffler modules. We will employ properties (P1)-(P3) to construct by induction on $\alpha$ a $\mathcal C _\kappa$-filtration $( M_\alpha \mid \alpha \leq \sigma )$ of $M$. This will show that $M \in \mathcal F \mathcal M _\kappa$. 

First, $M_0 = 0$. Assume that $M_\alpha \in \mathcal F \mathcal M _\kappa$ is defined so that $M_\alpha \otimes_R S \in \mathcal H$. Then either $M_\alpha = M$ and we let $\sigma = \alpha$, or else there is a $\leq \kappa$-generated submodule $X_0$ of $M$ such that $X_0 \nsubseteq M_\alpha$. By property (P3), there exists $H_1 \in \mathcal H$ such that $H_0 = M_\alpha \otimes_R S \subseteq (M_\alpha + X_0) \otimes_R S \subseteq H_1$ and $H_1/H_0$ is $\leq \kappa$-presented. So there exists a $\leq \kappa$ generated submodule $X_1$ of $M$ such that $X_0 \subseteq X_1$ and $H_1 \subseteq (M_\alpha + X_1) \otimes_R S$. Proceeding similarly, we obtain a chain of $\leq \kappa$-generated submodules $X_0 \subseteq X_1 \subseteq \dots$ of $M$ and a chain $H_0 \subseteq H_1 \subseteq \dots$ of elements of $\mathcal H $ such that 
$$H_0 = M_\alpha \otimes_R S \subseteq (M_\alpha + X_0) \otimes_R S \subseteq H_1 \subseteq (M_\alpha + X_1) \otimes_R S \subseteq H_2 \subseteq \dots$$
Let $X = \bigcup_{ < \omega} X_i$, $H = \bigcup_{i < \omega} H_i$, and $M_{\alpha +1} = M_\alpha + X$. Then $X$ is a $\leq \kappa$-generated submodule of $M$, so $M_{\alpha + 1}/M_\alpha$ is $\leq \kappa$-generated, too. Moreover, $H \in \mathcal H$ by property (P1), whence $M_{\alpha + 1} \otimes_R S = \bigcup_{i < \omega} (M_\alpha + X_i) \otimes_R S = H \in \mathcal H$. Since $\kappa$ is infinite and for each $i < \omega$, $H_{i+1}/H_i$ is a $\leq \kappa$-presented flat Mittag-Leffler $S$-module by property (P2), so is $H/H_0 = (M_{\alpha + 1} \otimes_R S)/(M_{\alpha} \otimes_R S)$. By Lemma \ref{descfML}, $M_{\alpha + 1}/M_\alpha$ is a $\leq \kappa$-generated flat Mittag-Leffler module, whence $M_{\alpha + 1}/M_\alpha \in \mathcal C _\kappa$ by Theorem \ref{absolut}(3). Thus $M_{\alpha +1} \in \mathcal F \mathcal M _\kappa$ and $M_{\alpha +1} \otimes_R S \in \mathcal H$.

If $\alpha$ is a limit ordinal, we let $M_\alpha = \bigcup_{\beta < \alpha} M_\beta$. Since the chain $(M _\beta \mid \beta < \alpha)$ is continuous, $M_{\alpha} \in \mathcal F \mathcal M _\kappa$. Moreover, $M_{\alpha} \otimes_R S = \bigcup_{\beta < \alpha} M_\beta \otimes_R S \in \mathcal H$ by property (P2). 
\end{proof}

\subsection{Finite type and the ad-property for flat relative Mittag-Leffler modules}\label{ftrel}

Now we turn to the setting of quasi-coherent sheaves arising from flat $\mathcal Q$-Mittag-Leffler modules. 

First, let us consider the boundary cases of $\mathcal Q = \{ 0 \}$ and $\mathcal Q = \lmod R$. As $\mathcal D _{\lmod R} = \mathcal F \mathcal M$, and being a flat Mittag-Leffler module is an ad-property by Lemma \ref{descfML}, the induced notion of a quasi-coherent sheaf is Zariski local for any scheme. The same holds for $\mathcal D _{\{0 \}} = \mathcal F _0$. 

In general, in order for the ad-property to hold, there has to be compatibility among the classes of left $R$-modules $\mathcal Q$ defining the meaning of {\lq}relative{\rq} for various rings $R$. The following results from \cite[\S 3]{BT2} will help us see what is needed:

\begin{lemma}\label{BT2rel} 
\begin{enumerate} 
\item Let $\varphi : R \to S$ be a flat ring homomorphism, $\mathcal Q \subseteq \lmod R$, and $M$ be a flat $\mathcal Q$-Mittag-Leffler module. Then $M \otimes_R S$ is a flat $\mathcal Q \otimes_R S$-Mittag-Leffler $S$-module.
\item Let $\varphi : R \to S$ be a faithfully flat ring homomorphism, and let $R \in \mathcal Q \subseteq \lmod R$. Assume that the implication 
$$(\star) \quad M \otimes_R S \hbox{ is a } (\mathcal Q \otimes_R S)\hbox{-Mittag-Leffler $S$-module} \implies M \in \mathcal D _{\mathcal Q}$$
holds for each countably presented flat module $M$. Then ($\star$) holds for each flat module $M$. 
\item Let $\mathcal S$ be a class of FP$_2$-modules and $\mathcal Q = \mathcal S ^\intercal$. Then the implication ($\star$) holds for every flat module $M$ and each faithfully flat ring homomorphism $\varphi : R \to S$. Moreover, for each flat ring homomorphism $\psi : R \to S$, we have $\Def {(\mathcal Q \otimes_R S)} = (\mathcal S \otimes_R S)^\perp$.    
\end{enumerate}
\end{lemma}   

Lemma \ref{BT2rel}(iii) suggests that one should concentrate on the case when $\mathcal Q$ are the classes of left $R$-modules of finite type from Example \ref{fintype}. The compatibility conditions sufficient for the ad-property are then as follows \cite[4.4]{BT2}: 

\begin{theorem}\label{mainBT2} For each ring $R$, let $\mathcal S _R$ be a class of FP$_2$-modules. Assume that the following compatibility conditions hold:
\begin{itemize}
\item[(C1)] $\mathcal S_R \otimes_R S \subseteq \mathcal S _S$ for each flat ring homomorphism $\varphi : R \to S$.
\item[(C2)] $(\mathcal S_R \otimes_R S )^\intercal = (\mathcal S _S)^\intercal$ for each faithfully flat ring homomorphism $\varphi : R \to S$.
\end{itemize}
Let $P_R$ be the property of being a flat $\mathcal Q _R$-Mittag-Leffler module, where $\mathcal Q _R = (\mathcal S _R)^\perp$. 

Then $P$ is an ad-property. In particular, the notion of a locally $P$-quasi-coherent sheaf is Zariski local.
\end{theorem}
\begin{proof} First, we prove the ascent along flat ring homomorphisms $\varphi : R \to S$. If $M \in \rmod R$ is flat and $\mathcal Q _R$-Mittag-Leffler, then $M \otimes_R S$ is a flat $(\mathcal Q \otimes_R S)$-Mittag-Leffler $S$-module by Lemma \ref{BT2rel}(i). By Proposition \ref{defenough} and Lemma \ref{BT2rel}(iii), $M \otimes_R S$ is a flat $(\mathcal S \otimes_R S)^\intercal$-Mittag-Leffler $S$-module. By Condition (C1), $(\mathcal S _S)^\intercal \subseteq (\mathcal S \otimes_R S)^\intercal$, whence $M \otimes_R S$ is a flat $\mathcal Q _S$-Mittag-Leffler $S$-module.

Let $\varphi : R \to S$ be a faithfully flat ring homomorphism and $M \in \rmod R$ be such that $M \otimes_R S$ is a flat $\mathcal Q _S$-Mittag-Leffler $S$-module. By Condition (C2) and Lemma \ref{BT2rel}(iii), $M \otimes_R S$ is a flat $\Def {(\mathcal Q \otimes_R S)}$-Mittag-Leffler $S$-module. Then $M$ is flat by (the proof of) Lemma \ref{descfML}, and $M$ is a $\mathcal Q _R$-Mittag-Leffler module by Lemma \ref{BT2rel}(iii).

The final assertion follows by Lemma \ref{acl}.
\end{proof}

The following is an application of Theorem \ref{mainBT2} to the particular case of locally f-projective quasi-coherent sheaves on \emph{coherent} schemes (i.e., those schemes whose structure sheaf consists of coherent rings):    

\begin{example}\label{coherentcont} We recall the setting of Example \ref{coherent}. So $\mathcal Q  = \{ R \}$, and the modules in $\mathcal D _{\mathcal Q}$ are called f-projective. 

Assume that $R$ is a coherent ring. Then $\Def {(\mathcal Q )}$ is the class of all flat modules, so $\Def {(\mathcal Q )} = (\mathcal S _R)^\intercal$ where $\mathcal S _R$ denotes the class of all finitely presented modules. Then condition (C1) clearly holds for each flat ring homomorphism of coherent rings.  

As for Condition (C2) in the setting of coherent rings, it suffices to prove that $(\mathcal S_R \otimes_R S )^\intercal \subseteq (\mathcal S _S)^\intercal$ for each faithfully flat ring homomorphism $\varphi : R \to S$. However, if $M \in (\mathcal S_R \otimes_R S )^\intercal$, then also $\Tor 1R{\mathcal S _R}M = 0$ by \cite[VI.4.1.1]{CE}, whence $M$ is a flat module, and $M \otimes_R S$ a flat $S$-module. Moreover, defining $f \in \Hom SM{M\otimes_R S}$ and $g \in \Hom S{M\otimes_RS}M$ by $f(m) = m \otimes 1$ and $g(m \otimes 1) = m$, we see that $gf = 1_M$, whence $M$ is isomorphic to a direct summand in $M \otimes_R S$. So $M$ is a flat $S$-module, and $M \in (\mathcal S _S)^\intercal$.    

Thus Theorem \ref{mainBT2} implies that the notion of a locally f-projective quasi-coherent sheaf is Zariski local for all coherent schemes.
\end{example}
  
\subsection{Locally tilting quasi-coherent sheaves}\label{zartilt}

Finally, we turn to Zariski locality in the settings induced by tilting modules. The results of this section come from \cite{HST}, and rely on the structure theory of tilting classes over commutative rings developed in \cite{APST}, \cite{Hr}, and \cite{HS}. 

Let $n \geq 0$. Consider the property $P_R$ of being an $n$-tilting $R$-module (see Definition \ref{tilting}). Thus, a \emph{locally $n$-tilting} quasi-coherent sheaf on a scheme $X$ is a quasi-coherent sheaf $\mathcal M$ such that $M(U)$ is an $n$-tilting $R(U)$-module for each open affine subset $U$ of $X$. 

A key tool for proving the Zariski locality in the tilting case is the following lemma from \cite[\S 2]{HST} which relies substantially on Theorem \ref{tilt}, that is, on tilting classes being of finite type. 
   
\begin{lemma}\label{primes} Let $\varphi : R \to S$ be a flat ring homomorphism, $T$ be an $n$-tilting module, $\mathcal A _T$ the induced left $n$-tilting class, $\mathcal B _T$ the induced $n$-tilting class, and $\mathcal S _T$ the representative set of all FP$_2$-modules in $\mathcal A _T$ (so that $\mathcal B _T = \mathcal S _T ^\perp$).

Let $T^\prime = T \otimes_R S$. Then $T^\prime$ is an $n$-tilting $S$-module, and $\mathcal B ^\prime = (\mathcal S _T \otimes_R S)^\perp$ is the $n$-tilting class induced by $T^\prime$. 
Let $\mathcal A^\prime = {}^\perp \mathcal B ^\prime$ be the left $n$-tilting class induced by $T^\prime$. 
\begin{itemize}
\item[(1)] $\mathcal A _T \otimes_R S \subseteq \mathcal A ^\prime$. Moreover, if $\varphi$ is faithfully flat, then for each module $M \in \rmod R$, $M \in \mathcal A _T$, iff $M \otimes_R S \in \mathcal A ^\prime$.
\item[(2)] $\mathcal B _T \otimes_R S \subseteq \mathcal B ^\prime$. Moreover, if $\varphi$ is faithfully flat, then for each module $M \in \rmod R$, $M \in \mathcal B _T$, iff $M \otimes_R S \in \mathcal B ^\prime$. 
\end{itemize}
\end{lemma}

The following lemma was proved in \cite[3.16]{HST}:    

\begin{lemma} \label{ffchange} Let $\varphi : R \to S$ be a faithfully flat ring homomorphism. Let $T^\prime$ be any $n$-tilting $S$-module of the form $M \otimes _R S$ for a module $M \in \rmod R$. Then there is an $n$-tilting module $T \in \rmod R$ such that $T \otimes _R S$ is equivalent to $T^\prime$. 
\end{lemma}

Now, we can prove our first claim concerning Zariski locality:

\begin{theorem}\label{adtilt} Let $n \geq 0$. 
\begin{itemize}
\item[(1)] The property of being an $n$-tilting module ascends along flat ring homomorphisms. 
\item[(2)] If $\varphi : R \to S$ is a faithfully flat ring homomorphism, and $\bar T$ is a module such that $\bar{T}^\prime = \bar T \otimes _R S$ is an $n$-tilting $S$-module, then $\bar T$ satisfies conditions (T1) and (T2).   
\end{itemize}
\end{theorem}
\begin{proof} (1) Let $\varphi : R \to S$ be a flat ring homomorphism and $T$ be an $n$-tilting module. Then $T \otimes _R S$ is an $n$-tilting $S$-module by Lemma \ref{primes}.

(2) Let $\varphi : R \to S$ be a faithfully flat ring homomorphism and $\bar T$ be a module such that $\bar{T}^\prime = \bar T \otimes _R S$ is an $n$-tilting $S$-module. By Lemma \ref{ffchange}, there is an $n$-tilting module $T$ such that $\bar{T}^\prime$ is equivalent to the $n$-tilting $S$-module $T ^\prime = T \otimes _R S$. Let $\mathcal A ^\prime$ and $\mathcal B ^\prime$ be the $n$-tilting classes induced by $T^\prime$ (equivalently, by $\bar{T}^\prime$) in $\rmod S$. Then $\bar T \otimes _R S \in \Add {T^\prime} = \mathcal A ^\prime \cap \mathcal B ^\prime$.

By Lemma \ref{primes}, $\bar T \in \mathcal A _T \cap \mathcal B _T = \Add T$. Since conditions (T1) and (T2) hold true for $T$, they also hold for $\bar T$.
\end{proof}

It is an open problem whether the property of being an $n$-tilting module descends along all faithfully flat ring homomorphisms. In \cite{HST}, a positive answer was given for $n \leq 1$, and for the case of faithfully flat ring homomorphisms of commutative \emph{noetherian} rings. 

However, Zariski locality does hold in general, as one only needs to prove the descent along the particular faithfully flat ring homomorphisms of the form $\varphi : R \to S = \prod_{j < m} R[f_j^{-1}]$ where $R = \sum_{j<m} f_jR$ (see \ref{acl}). Moreover, in the presence of conditions (T1) and (T2), condition (T3) can be replaced by a homological condition involving the unbounded derived category $D(R)$ of $R$-modules. Namely, (T3) is then equivalent to the condition that $D(R)$ is the smallest localizing subcategory of itself containing $T$. The latter condition can be verified in the setting of Theorem \ref{adtilt}(2) (see \cite[Lemma 4.1]{HST} for more details). Thus we conclude:

\begin{theorem}\label{final} Let $n \geq 0$. Then the notion of a locally $n$-tilting quasi-coherent sheaf is Zariski local for all schemes. 
\end{theorem}
       
\begin{remark}\label{correct}
Unlike the previous sections, the case of $n = 0$ in Theorem \ref{final} differs from the case of vector bundles. Namely, $0$-tilting modules are exactly the (possibly infinitely generated) projective \emph{generators}. However, Lemmas \ref{primes} and \ref{ffchange} for $n = 0$ do imply that the property of being a projective module is an ad-property. The point is that the assumption that any $n$-tilting $S$-module $T^\prime$ is equivalent to an $n$-tilting $S$-module of the form $M \otimes _R S$ for a module $M \in \rmod R$ is satisfied for $n = 0$. 

However, this assumption fails for $n \geq 1$: by \cite[6.2]{HS}, $n$-tilting classes in $\rmod R$ correspond 1-1 to certain $n$-tuples of subsets of $\Spec R$ called \emph{characteristic sequences}. The existence of a faithfully flat ring homomorphism $\varphi : R \to S$ only gives a monomorphism from the characteristic sequences in $\Spec R$ to those in $\Spec S$, not a bijection (see \cite[3.8]{HST} for more details).     
\end{remark}

\end{document}